\documentclass[letterpaper,10pt,3p,preprint]{elsarticle}

\usepackage{amsmath,amssymb,amsfonts}
\usepackage{bm}
\usepackage{tikz-cd}
\usepackage{graphicx}

\usepackage[hidelinks]{hyperref}

\newcommand{\Nbb}{\mathbb{N}}
\newcommand{\Rbb}{\mathbb{R}}
\newcommand{\Cbb}{\mathbb{C}}
\newcommand{\Ebb}{\mathbb{E}}

\newcommand{\Fcal}{\mathcal{F}}

\DeclareMathOperator{\ran}{ran}
\DeclareMathOperator{\tr}{tr}
\DeclareMathOperator{\Id}{Id}
\DeclareMathOperator{\proj}{proj}

\newcommand{\vect}[1]{\bm{#1}}
\newcommand{\matr}[1]{\bm{#1}}

\begin{document}

\begin{frontmatter}

\title{Quantum mechanical closure of partial differential
equations with symmetries}

\author[1]{Chris Vales\corref{cor1}}
\ead{chris.vales@dartmouth.edu}
\author[1]{David C. Freeman}
\ead{david.c.freeman.gr@dartmouth.edu}
\author[1]{Joanna Slawinska}
\ead{joanna.m.slawinska@dartmouth.edu}
\author[1]{Dimitrios Giannakis}
\ead{dimitrios.giannakis@dartmouth.edu}

\cortext[cor1]{Corresponding author.}
\affiliation[1]{organization={Department of Mathematics,
    Dartmouth College},
    city={Hanover},
    state={NH 03755},
    country={USA}}

\begin{abstract}
We develop a statistical framework for the dynamical closure of
spatiotemporal dynamics governed by partial differential equations.
Employing the mathematical framework of quantum
mechanics to embed the original classical dynamics into a quantum
mechanical representation, we use the space of quantum density
operators to model the unresolved degrees of freedom of the
original dynamics in a statistical sense, and the framework
of quantum measurement to predict their contributions to the
resolved dynamics.
The embedded dynamics is discretized by a positivity preserving
process, leading to a compressed representation that is invariant
under the dynamical symmetries of the resolved dynamics.
We present a data based formulation of the closure scheme
and apply it to a closure problem for the shallow water equations.
The numerical results demonstrate that our closure model can accurately
predict the main features of the true dynamics, including for out
of sample initial conditions.
\end{abstract}


\end{frontmatter}

\section{Introduction}\label{sec:intro}
Complex dynamical systems often encompass degrees of freedom that
evolve over a wide range of temporal and spatial scales, with
their evolution potentially depending on different physical regimes.
A challenging problem in such multiscale, multiphysics systems
is incorporating dynamical information for degrees of freedom
that either are too computationally expensive to simulate directly,
or depend on fully or partially unknown physical processes
\cite{Sagaut2006,WE2011,Stensrud2013}.

The approach employed by \emph{closure} or \emph{parametrization}
methods is to derive surrogate dynamical models that can be used
to approximate the contributions of fine grain degrees of freedom
to the dynamics of coarse grain processes we can simulate
consistently.
Early applications included turbulent flows and climate modeling
\cite{Arakawa1974,Gent1990,Sagaut2006},
but closure methods have since found consequential applications
in a wide range of fields \cite{WE2007}.
For example, in atmospheric dynamics a technique known as
super-parametrization employs column models with explicitly resolved
cloud microphysics as surrogate models of small scale moist convection,
embedded within the discretization cells of a coarse atmospheric model
\cite{Grabowski1999,Grooms2013}.
In plasma dynamics there are ongoing efforts to derive models
for subgrid scale processes in fluid models,
as well as to derive fluid closures of kinetic models,
incorporating both physical modeling and data based approaches
\cite{Hammett1990,Joseph2016,ChMa2020,Laperre2022}.

In this work we focus on the development of a statistical,
data based framework for the closure of dynamics governed by
partial differential equations (PDEs).
More specifically, we employ the operator theoretic framework of
quantum mechanics to define surrogate dynamical models, using quantum
density operators to encode statistical information about fine grain
degrees of freedom and quantum observables to model their contributions
to the coarse dynamics.
Our approach extends a previously developed quantum mechanical
closure (QMCl) framework \cite{Freeman2024}
to the closure of spatiotemporal dynamics, explicitly taking into
account the presence of dynamical symmetries to improve the
numerical efficiency and training data requirements of our
surrogate models \cite{Giannakis2019vsa}.

In the rest of Section \ref{sec:intro}
we briefly discuss the problem of dynamical closure and
different modeling approaches toward addressing it,
present the main aspects of our proposed closure scheme
and summarize the main results of this work.
In Section \ref{sec:qmcl} we develop in detail the different
aspects of our closure framework.
In Section \ref{sec:swe} we present an application to a closure
problem for the shallow water equations,
followed by a brief conclusion in
Section \ref{sec:conclusion}.
Finally, \ref{app:implementation} provides additional details
on the numerical implementation of our closure scheme.
The code used to generate the numerical results and plots of this
work is available
online\footnote{\url{https://github.com/cval26/qmcl_swe}}.

\subsection{Dynamical closure}\label{sec:intro-closure}
We consider a dynamical system $\Phi^t\colon X\rightarrow X$,
$t\geq 0$ on a state space $X$.
The flow map $\Phi^t$ represents the spatiotemporal dynamics
governed by a PDE on a bounded spatial domain $S$.
As such, the state space $X$ is a Hilbert space of functions
on the considered spatial domain.
To define the general dynamical closure problem, we decompose the
state space as $X=X_r\oplus X_u$,
with $X_r$ denoting the state space of the resolved degrees
of freedom of the dynamics and $X_u$ the state space of the
unresolved degrees.
The flow map of the original dynamics $\Phi^t$
can then also be decomposed into the resolved dynamics
$\Phi_r^t\colon X\to X_r$
and the unresolved dynamics
$\Phi_u^t\colon X\to X_u$.

To provide an illustrative example of this decomposition we consider
the advection equation with periodic boundary conditions
\begin{equation}\label{eq:model-pde}
    \partial_tu+c\partial_xu=F(u)
\end{equation}
where $u$ denotes a function of time $t\geq 0$ and
space $x\in S\subset\Rbb$, $F$ a nonlinear function of $u$,
and $c\in\Rbb$ the advection velocity.
We define the spatially averaged state function
\begin{equation*}
    \hat{u}(\cdot,x)=\frac{1}{\ell}\int_{x-\ell/2}^{x+\ell/2}u(\cdot,y)dy
\end{equation*}
where $\ell>0$ denotes the length of the averaging window.
The evolution equation for $\hat{u}$ is derived by spatially averaging
\eqref{eq:model-pde} to obtain
\begin{equation}\label{eq:model-pde-coarse}
    \partial_t\hat{u}+c\partial_x\hat{u}=F(\hat{u})+G(\hat{u},u)
\end{equation}
with residual term
\begin{equation}\label{eq:model-pde-residual}
    G(\hat{u},u)(\cdot,x)=\frac{1}{\ell}\int_{x-\ell/2}^{x+\ell/2}
        F(u)(\cdot,y)dy-F(\hat{u})(\cdot,x).
\end{equation}
A version of this type of spatial averaging on a discretized spatial
domain will be used for our numerical application in
Section \ref{sec:swe}.

We consider the original PDE \eqref{eq:model-pde}
as the equation generating the true dynamics $\Phi^t$
with $u\in X$ as the state variable.
Motivating the decomposition described above, we take the
spatially averaged equation \eqref{eq:model-pde-coarse}
as the equation governing the resolved dynamics, 
meaning that $\hat{u}$ is the resolved state variable with
state space $X_r$ and $\Phi_r^t$ the flow map governing its evolution.
The associated unresolved state space $X_u$ is then implicitly defined
as the orthogonal complement of $X_r$ in $X$,
representing phenomena that cannot be resolved by $\hat{u}$.

Due to the nonlinear term $F(u)$ in the original PDE, the
spatial averaging operation has not eliminated the dependence of
\eqref{eq:model-pde-coarse}
on the unresolved degrees of freedom.
This is clearly demonstrated by the presence of the residual terms
\eqref{eq:model-pde-residual},
which retain their dependence on the original state variable $u$.
As a result, the resolved dynamics map $\Phi_r$ generated by
\eqref{eq:model-pde-coarse} requires at least partial knowledge of
unresolved variables at any given time.
This is a manifestation of the \emph{dynamical closure problem}:
to form a closed system of evolution for the resolved degrees of freedom,
we need knowledge of additional, unresolved variables.
The residual terms \eqref{eq:model-pde-residual} still depending on the
unresolved degrees of freedom are often referred to as the flux terms
or \emph{fluxes}.

To address this problem, \emph{closure} or \emph{parametrization methods}
construct a surrogate model for the resolved degrees of freedom,
which approximates the evaluation of the fluxes by modeling the
unresolved variables appropriately.
The surrogate models formed for the resolved and unresolved variables
can then evolve in parallel, leading to a closed system of evolution
for the resolved variables.
To be useful, the model employed for the unresolved variables has to
be computationally efficient and lead to surrogate dynamics that
constitutes an effective and consistent approximation of the original
dynamics in an appropriate sense.

Closure problems in physics also arise when deriving a fluid model
from a kinetic description of a system by considering the evolution of
moments of the phase space probability density function
\cite{Bardos1991}.
In that setting, the resolved variables are the low order velocity
moments involved in the fluid description---usually mass, momentum
and energy---while the unresolved variables are the ignored
higher order moments.
To define a closed system of evolution for the fluid model,
one is again tasked with forming a surrogate model to estimate the
required fluxes.

\subsection{Modeling approaches}\label{sec:intro-modeling}
One approach toward defining a surrogate model is to model the
unresolved variables by a deterministic function of the resolved
variables.
Defining such a function can be achieved by both physical modeling
and data based approaches
\cite{Stensrud2013,Rasp2018,Yuval2020,Bae2022,Bracco2025}.
A drawback of this approach is that it corresponds to a reduction
in the dimension of the full state space of the dynamics, thereby
potentially reducing the overall complexity that can be produced
by the surrogate model \cite{Palmer2001}.

Our work is motivated by an alternative, statistical approach,
where the state space of the unresolved degrees of freedom is
replaced by the space of probability density functions on
the full state space of the dynamics
\cite{Majda1999,Chorin2015,Berner2017}.
More specifically, a probability density function is used to assign
a probability that the dynamical system is in a given state at
a given time.
Contrary to purely functional approaches, where a single state of the
system is uniquely identified for any given time, statistical methods
take into consideration the presence of uncertainty in our knowledge
of the full state of the dynamical system.

To illustrate the main features of this approach we define a probability
density function $p\in L^1(X)$ on the state space.
This is a real valued, nonnegative function whose integral over $X$
is equal to one with respect to an appropriate measure.
At any given time, the evaluation of the density function $p(u)\geq 0$
assigns a probability that the dynamical system is currently in state
$u\in X$.
Equipped with statistical knowledge of the current state of the dynamics,
one can use it to compute the fluxes required for the evolution of
the resolved variables.
To that end, the fluxes are modeled by a bounded, real valued function
$f\in L^\infty(X)$ called an observable; namely, given a state
$u\in X$, the associated flux is $f(u)\in\Rbb$.
The required flux can then be estimated by the expected value of $f$
with respect to $p$, namely $\int_Xf(u)p(u)du\in\Rbb$ with respect
to an appropriate measure.

Having estimated the required flux term, one can then proceed with the
evolution of the resolved variables; in our motivating example, this
corresponds to integrating in time the spatially averaged PDE.
For the surrogate evolution system to be closed,
it must also account for the time evolution of the surrogate model
of the unresolved variables, namely the density function $p$.
This can be accomplished by using the transfer or Perron-Frobenius
operator, which governs the evolution of probability density functions
under the dynamics \cite{Lasota1994}.
These notions will be made precise in the following sections.

As explained in the following section, one of the main features
differentiating our closure framework from the general statistical
framework sketched above is that we replace the probability density
function $p$ with a probability density operator.
Namely, the nonnegative function $p$ with unit integral is replaced by
a positive semidefinite operator with unit trace.
Accordingly, the bounded function $f$ modeling the fluxes is replaced
by a bounded operator.
In our framework, the estimated flux is then given by the expected
value of the observable operator with respect to the probability
density operator.

\subsection{Quantum mechanical closure framework}\label{sec:intro-qmcl}
We provide a brief overview of the QMCl framework as it has been
developed in \cite{Freeman2024}
for the closure of dynamics governed by ordinary differential equations.
A detailed treatment of the framework will be given in
Section \ref{sec:qmcl},
where we extend it to the closure of dynamics governed by PDEs.

As already mentioned earlier, we adopt a
statistical description of the unresolved dynamics, meaning that we
model the evolution of a probability density function of the state
vectors, rather than the evolution of the state vectors themselves.
One of the main features of the QMCl framework is that this statistical
description is embedded into the operator setting of quantum mechanics
on a complex Hilbert space $H=L^2(X,\mu)$, where $X$ denotes the
state space and $\mu$ an invariant probability measure of the
original dynamics.
This involves embedding observables and probability densities,
as well as the evolution of said densities under the dynamics and
their bayesian conditioning.

In the classical statistical description, a probability
density function $p\in L^1(X,\mu)$ was used
to model the state of the unresolved degrees of freedom.
Given that $\sqrt{p}\in H$, we now map function $p$ to the
quantum density operator $\rho\in B(H)$,
$\rho=\langle\sqrt{p},\cdot\rangle\sqrt{p}$,
representing a pure quantum state.
Quantum density operators---also called quantum states---are
positive semidefinite, selfadjoint operators in $B(H)$
with unit trace, and can be thought of as the operator analog
of probability density functions.
In this work we mainly use pure quantum states,
which are rank-one quantum density operators equivalent
to projections along a unit vector
\cite{Hall2013,Sakurai2021}.
The statistical description of the dynamics obtained by mapping the
density function $p$ to the density operator $\rho$ is an example
of the Koopman-von Neumann representation of classical dynamics
\cite{Mauro2002,Bondar2019,Joseph2020}.

The fluxes needed for the closure of the resolved dynamics were
modeled classically by an observable
$f\in L^\infty(X,\mu)\subset H$.
This is now mapped to a selfadjoint multiplication operator
$A\in B(H)$, $Ag=fg$ for all $g\in H$,
which is referred to as a quantum observable.
With these definitions in place, the surrogate flux needed for the
evolution of the resolved dynamics is defined as the expected value
of $A$ with respect to $\rho$, namely
$\Ebb_\rho A=\tr(\rho A)\in\Rbb$,
which is connected to the theory of quantum measurement.
In other words, we use the space of quantum density operators
to model the unresolved degrees of freedom of the dynamics
and the framework of quantum measurement to estimate their
fluxes.

In addition to evaluating the required flux terms, we need to be
able to evolve the density operator $\rho$ under the system's
dynamics and condition it based on observed dynamical behavior.
In the classical case, the square root of the probability density
function $\sqrt{p}$ evolved according to the transfer or
Perron-Frobenius operator $P\colon H\to H$,
$P\sqrt{p}=\sqrt{p}\circ\Phi^{-1}$.
Accordingly, the evolution of the density operator $\rho$
is governed by the induced transfer operator
$\mathcal{P}\colon B(H)\to B(H)$,
acting by conjugation $\mathcal{P}\rho=P\rho P^{-1}$.

In analogy to the bayesian conditioning of the probability density
function $p$, we can condition the density operator $\rho$
using a quantum effect operator $e\in B(H)$
\begin{equation}\label{eq:quantum-bayes}
    \rho\vert_e=\frac{\sqrt{e}\rho\sqrt{e}}{\Ebb_\rho e}
\end{equation}
where $\sqrt{e}$ denotes the positive square root of $e$.
The above equation can be interpreted as the operator
analog of the classical Bayes rule \cite{Busch1996}.
In addition, the quantum effect operator $e$ can be thought of as
the operator analog of an indicator function, which is used in
classical probability to represent an event, namely a set in the
event space.
For instance, for a given set $F\subset X$ we can build the
indicator function $1_F\in H$ taking the value one on $F$
and zero elsewhere.
The quantum effect operator $e_F\in B(H)$ can then be defined
as multiplication by $1_F$, namely $e_Fg=1_Fg$ for all $g\in H$.

The following diagram depicts the evolution for one timestep of
the QMCl closure model built using the tools developed in the
previous paragraphs.
The top row of the diagram is concerned with the resolved degrees
of freedom $\hat{u}\in X_r$, while the bottom row with the density
operator $\rho$ acting as our surrogate model for the unresolved
degrees of freedom.
In the diagram, subscripts are used to denote the time
index, while $\Id$ denotes the identity map.
\begin{equation}\label{dgm:intro-qmcl}
\begin{tikzcd}[column sep=huge, row sep=large]
    \hat{u}_n\arrow[r,"\text{resolved}"]
        &\hat{u}_{n+1}\arrow[r,"\Id"]\arrow[rd,"\text{effect}"]
        &\hat{u}_{n+1}\\
    \rho_n\arrow[ru,"\text{flux}"]\arrow[r,"\text{transfer}"']
        &\rho_{n+1}'\arrow[r,"\text{Bayes}"'] &\rho_{n+1}
\end{tikzcd}
\end{equation}

To evolve the resolved variables for one timestep, we require
knowledge of the flux term depending on the unresolved degrees
of freedom.
As explained above, this is estimated by the expected value
of an observable $A$ with respect to the density $\rho_n$.
The evolution of the density operator $\rho_n$ consists of two
steps, a prediction step and a correction step.
The prediction step employs the transfer operator to update the
density operator $\rho_n$, yielding the prior density
$\rho_{n+1}'$.
The correction step relies on the quantum Bayes formula
$\eqref{eq:quantum-bayes}$ to correct this prediction,
yielding the posterior density $\rho_{n+1}$.
The bayesian correction uses information extracted from the updated
resolved variables $\hat{u}_{n+1}$, in the form of a quantum effect,
to condition the prior density, concluding one timestep of the
evolution of the density operator.

\subsection{Main results}
We offer a preliminary description of the way QMCl is extended
to the closure of dynamics governed by PDEs.
As explained in Section \ref{sec:intro-closure},
in the setting of spatiotemporal dynamics the resolved state
$\hat{u}$ is a spatially averaged function on the spatial domain $S$.
Motivated by this, the main point of departure of the extended
QMCl framework from the one outlined above is that we use a field
of density operators $\rho(x)$ to model the unresolved degrees
of freedom at each point $x$, instead of using one density operator
for the whole domain.
This enables us to compute surrogate fluxes that vary over the spatial
domain, and are thus capable of modeling the true flux terms
\eqref{eq:model-pde-residual} introduced earlier.
Schematically, this means that the extended QMCl scheme involves
one copy of diagram \eqref{dgm:intro-qmcl} for each gridpoint
$x$, with the diagrams evolving in parallel.

The main features of the present work can be summarized as follows.
\begin{enumerate}
\item\emph{Positivity preservation}.
    Using probability density operators to encode statistical
    information about the original dynamics leads to a closure
    scheme that is positivity preserving, which is important for
    sign definite quantities.
    More specifically, positive flux terms in the original dynamics are
    represented by positive definite operators, in turn leading to
    positive surrogate flux terms in our scheme's implementation
    (Section \ref{sec:qmcl-data}).
\item\emph{Symmetry factorization}.
    We incorporate elements of the vector valued spectral analysis (VSA)
    framework \cite{Giannakis2019vsa}
    to develop a data based formulation of our scheme that is efficient
    and effective in the presence of dynamical symmetries.
    More specifically, our scheme rests on a compressed approximate
    representation of the original dynamics that is invariant with
    respect to the dynamical symmetries of the resolved dynamics
    (Section \ref{sec:qmcl-symmetries}).
\item\emph{Computational efficiency}.
    The numerical implementation of our closure scheme
    requires solving an eigenvalue problem for a kernel matrix,
    which is one of the most expensive tasks of the scheme
    (Section \ref{sec:qmcl-data}).
    To address that we compute a low rank approximation of the original
    kernel matrix and use it to solve the eigenvalue problem with reduced
    cost \cite{Vales2025evd}.
\item\emph{Application}.
    We apply our scheme to a closure problem for the shallow water
    equations on a one dimensional spatial domain
    (Section \ref{sec:swe}).
    The results demonstrate that the derived surrogate model can
    effectively capture and reproduce the main features of the dynamics,
    including for out of sample initial conditions.
    The numerical results can be reproduced using the provided
    code\footnote{\url{https://github.com/cval26/qmcl_swe}}.
\end{enumerate}

We end this section by offering a partial answer to the question:
why formulate a statistical closure framework using the operator
setting of quantum mechanics, instead of the original function setting
of classical mechanics?
First, working in an operator setting leads to an implementation
of our closure scheme that is naturally positivity preserving, meaning
that the required discretization of classical observables and
probability density functions preserves their positivity
(Section \ref{sec:qmcl-data}).
Second, the use of density operators to encode statistical information
about the unresolved degrees of freedom offers a potential increase in
the available encoding capacity, compared to using density functions
(Section \ref{sec:swe-discussion}).
However, these benefits are accompanied by an increase in computation
and memory cost, since most discretized quantities are now represented
by matrices instead of vectors
(Section \ref{sec:qmcl-data}).
Finally, the fact that our closure scheme rests on the mathematical
framework of quantum mechanics offers the opportunity to explore
its future implementation in quantum computers, especially in
order to investigate whether a more efficient implementation can be
devised by doing so.

\section{Closure of partial differential equations}\label{sec:qmcl}
In this section we develop the extended QMCl closure scheme and
its data based formulation.
We begin by setting up the problem in
Section \ref{sec:qmcl-setting},
before presenting the main aspects of our closure modeling
framework in Section \ref{sec:qmcl-infinite}.
As explained in Section \ref{sec:qmcl-infinite},
the extended QMCl scheme is formulated in an infinite
dimensional Hilbert space $H$.
In Section \ref{sec:qmcl-finite} we construct
a finite dimensional subspace of $H$, which will be used
in the numerical implementation of our scheme.
Section \ref{sec:qmcl-symmetries} then investigates some of
the properties of the constructed subspace with respect to
dynamical symmetries.
Finally, Sections \ref{sec:qmcl-data} and \ref{sec:qmcl-multi}
deal with the data based formulation of the extended QMCl scheme
that will be used for our numerical experiments.

\subsection{Problem setting}\label{sec:qmcl-setting}
Expanding on the preliminary material presented in
Section \ref{sec:intro-closure},
we consider the spatiotemporal dynamics
$\Phi^n\colon X\to X$, $n\in\Nbb$,
possessing an invariant probability measure $\mu$.
The discrete time map $\Phi^n$ is assumed to be a discrete time
version of the continuous time dynamics, with
$\Phi^n=\Phi^{n\Delta t}$ for a timestep $\Delta t>0$.
Based on the guiding example of the advection equation given in
Section \ref{sec:intro-closure},
the dynamics map $\Phi^n$ is generated by the temporally discretized
version of \eqref{eq:model-pde-coarse}, corresponding to the
spatially averaged dynamics.

The state space $X$ of the dynamics is a Hilbert space of functions
on a bounded spatial domain $S\subset\Rbb^c$, $c\in\Nbb$.
Specifically, we assume
$X\subset L^2(S,\nu;\Rbb^d)$, $d\in\Nbb$,
with $\nu$ denoting the Lebesgue measure.
We consider the state space decomposition $X=X_r\oplus X_u$
into the resolved and unresolved state spaces.
Based on \eqref{eq:model-pde}--\eqref{eq:model-pde-coarse},
the resolved state space is associated with spatially averaged
or coarsened functions, which were denoted by $\hat{u}$ in
Section \ref{sec:intro-closure}.
Analogously, the unresolved state space $X_u=X_r^\perp$ is defined
as the orthogonal complement, representing degrees of freedom that
cannot be resolved by the coarse state $\hat{u}$.

Introducing the product state space $\Omega=X\times S$
and product measure $\sigma=\mu\times\nu$,
we define the Hilbert space of spatiotemporal observables
$H=L^2(\Omega,\sigma;\Cbb)$.
Given a real valued observable $f\in H$, for almost every
$u\in X$ we have that $f(u,\cdot)\in L^2(S,\nu;\Rbb)$
yields a real valued function on the spatial domain $S$,
with $f(u,x)\in\Rbb$ its value at $x\in S$.
The motivation for our choice of space $H$ is our interest
in spatiotemporal observables; that is, patterns of evolution
that have both a temporal dependence through $u\in X$
and a spatial dependence through $x\in S$.
An example of such an observable is the flux function introduced
in \eqref{eq:model-pde-residual};
given a state $u\in X$ and gridpoint $x\in S$,
it outputs the residual term required for the local closure
of the coarse dynamics \eqref{eq:model-pde-coarse}.
Moving forward, the range space of a function space is assumed
to be $\Cbb$ unless explicitly noted otherwise, and the measure
notation is dropped when clear from the context; for example,
we write $L^2(X)$ instead of $L^2(X,\mu;\Cbb)$.

The fact that we employ an invariant measure $\mu$ means
that the dynamics $\Phi^n$ is volume preserving---also called
incompressible---when the state space ``volume'' is measured
with $\mu$.
This means that $\mu(\Phi^{-n}(F))=\mu(F)$
for every measurable set $F\subseteq X$ and $n\in\Nbb$.
For instance, there are dissipative nonlinear PDEs which, although
not volume preserving with respect to the ambient Lebesgue measure,
can be shown to possess an invariant measure related to their
dynamics on an attractor \cite{Temam1997,Robinson2001}.
The existence of an invariant probability measure for compact $X$
and continuous $\Phi^n$ is ensured by the Krylov-Bogoliubov
theorem \cite{Lasota1994},
which is sufficient for the applications considered in the
present work.

In the data based formulation of our scheme treated in
Sections \ref{sec:qmcl-data}, \ref{sec:qmcl-multi} and \ref{sec:swe},
we use dynamical trajectories to collect state space samples as
our training data.
The process of sampling the dynamics along a trajectory can be
naturally interpreted as sampling an invariant measure $\mu$,
since it is the dynamics that dictates the statistical distribution
of the sampled points in the state space.
This may be contrasted to sampling the dynamics on a predefined
cartesian grid in the state space, which would correspond to
sampling the ambient Lebesgue measure.

\subsection{Closure model}\label{sec:qmcl-infinite}
The flux terms needed for the closure of the resolved dynamics
are represented by $d$ classical real valued observables
$f_j\in L^\infty(\Omega)\subset H$, $j\in\{0,\ldots,d-1\}$.
Following the approach outlined in Section \ref{sec:intro-qmcl},
each observable $f_j$ is mapped to a bounded, selfadjoint
multiplication operator $A_j\in B(H)$, which we refer to as
a quantum observable.

As mentioned earlier, an example of such an observable is
the flux function \eqref{eq:model-pde-residual},
which generates the residual terms needed for the integration
of the coarse dynamics \eqref{eq:model-pde-coarse}.
Given that the advection equation \eqref{eq:model-pde}
models the evolution of the scalar field $u$, it follows that
$d=1$ and thus only one observable is needed to model the fluxes.
In Section \ref{sec:swe} we are going to consider a closure
problem for the shallow water equations \eqref{eq:swe},
a system of PDEs for two scalar fields (height and momentum).
As a result, in that case $d=2$ and thus two observables
will be used to generate the needed fluxes.

To model the unresolved degrees of freedom we employ a field
of quantum density operators $\rho\colon S\to B(H)$
over the discretized spatial domain $S$.
Each density operator $\rho(x)\in B(H)$ represents the unresolved
degrees of freedom of the dynamics at a point $x\in S$
in a statistical sense.
It is defined as the rank-one projection along the square root
of a probability density function $p_x\in L^1(\Omega)$,
$\rho(x)=\langle\sqrt{p_x},\cdot\rangle\sqrt{p_x}$.
Using a field of density operators allows us to compute fluxes that
vary over the spatial domain $S$ for each quantum observable $A_j$.
More specifically, we define the corresponding
$d$ surrogate flux functions
$\tilde{f}_j\colon S\to\Rbb$
with action
\begin{equation}\label{eq:surrogate-flux}
    \tilde{f}_j(x)=\tr(\rho(x)A_j)\in\Rbb
\end{equation}
yielding $d$ surrogate flux terms at each gridpoint.

For the closure of the coarse advection dynamics
\eqref{eq:model-pde-coarse}, one density operator $\rho(x)$
is required for each gridpoint $x$ of the spatial domain.
Using equation \eqref{eq:surrogate-flux}
with the single quantum observable modeling the flux terms
\eqref{eq:model-pde-residual},
we can generate one surrogate flux at each gridpoint $x$
for each timestep, thus creating a surrogate evolution equation
for the coarse dynamics.

The time evolution of each density operator $\rho(x)$
is performed by conjugation with the transfer operator
$P\colon H\to H$,
$P\sqrt{p_x}\,(u,\cdot)=\sqrt{p_x}\,(\Phi^{-1}(u),\cdot)$,
yielding $\rho^{n+1}(x)=P\rho^n(x)P^{-1}$,
where superscripts denote the time index.
Finally, the bayesian conditioning of each density operator
$\rho(x)$ is performed using a quantum effect operator
$e_x\in B(H)$ and the quantum Bayes rule
\eqref{eq:quantum-bayes} introduced in
Section \ref{sec:intro-qmcl}.
The way the quantum effects $e_x\in B(H)$ are computed in practice
will be explained in Section \ref{sec:qmcl-data}.

\subsection{Finite dimensional projection}\label{sec:qmcl-finite}
As explained above, our extended closure model is formed in the
infinite dimensional Hilbert space of observables $H$.
In this section we form a finite dimensional subspace
$H^L\subset H$ by constructing a collection of orthonormal basis
functions that span $H^L$.

To construct a finite dimensional subspace of $H$ we use the
orthonormal basis generated by the eigenfunctions of a
positive definite, selfadjoint kernel integral operator
$K\colon H\to H$
\begin{equation}\label{eq:integral-op}
    Kf=\int_\Omega \kappa(\cdot,\omega)f(\omega)d\sigma(\omega)
\end{equation}
with a positive definite, symmetric and uniformly continuous
kernel function
$\kappa\in C_b(\Omega\times\Omega)\cap
L^2(\Omega\times\Omega)$.
By forming an operator acting directly on the product Hilbert
space $H$ we obtain eigenfunctions that are generally not of
separable, tensor product form.
This is in contrast to traditional approaches such as the
proper orthogonal decomposition, where one would build a basis
for $L^2(X)$ (the temporal space in this context)
and form its tensor product with a basis for the
spatial space $L^2(S)$ to form a basis for
$H\simeq L^2(X)\otimes L^2(S)$.
As explained below, in the presence of dynamical symmetries
our approach leads to a more compressed finite dimensional
representation of the dynamics on $H$,
meaning that each eigenfunction of operator $K$
as defined above can represent spatiotemporal patterns that
are not expressible as a product of a single pair of temporal
and spatial modes \cite{Giannakis2019vsa}.

The kernel function
$\kappa(\omega,\omega')=k(W_Q(\omega),W_Q(\omega'))$
is formed by composing a gaussian kernel function
$k\colon\Rbb^{dQ}\times\Rbb^{dQ}\to\Rbb$
and a time delay embedding map
$W_Q:\Omega\to\Rbb^{dQ}$,
where $Q\in\Nbb$ denotes the number of time delays.
More specifically, the map $W_Q$ acts by
\begin{equation}\label{eq:delay-map}
    W_Q((u,x))=\bigl(\proj_r(u)(x),\,
        \proj_r(\Phi^{-1}(u))(x),\,\ldots,\,
        \proj_r(\Phi^{-(Q-1)}(u))(x)\bigr)
\end{equation}
which corresponds to pointwise time delay embedding of
the product state sample $\omega=(u,x)$.
Note that $W_Q$ considers only the resolved component
of state vector $u\in X$ by using the projection operator
$\proj_r\colon X\to X_r$.
This is because in many practical situations one has access to data
only for the resolved variables, not the full state of the system.
The gaussian kernel then acts by
\begin{equation}\label{eq:data-kernel}
    k(W_Q(\omega),W_Q(\omega'))=\exp\Bigl(-\frac{1}{\epsilon Q}
        \lVert W_Q(\omega)-W_Q(\omega')\rVert^2\Bigr)
\end{equation}
where $\lVert\cdot\rVert$ denotes the euclidean norm in $\Rbb^{dQ}$
and $\epsilon>0$ is a tunable bandwidth parameter.
The bandwidth parameter $\epsilon$ can be calibrated using the
procedure developed in \cite{Coifman2008}.
In our numerical experiments we normalize the above kernel to turn it
into a bistochastic kernel.
The bistochastic kernel normalization process and motivation
for using it can be found in
\cite{Coifman2013,Das2021}.

The eigenfunctions of operator $K$ from \eqref{eq:integral-op}
form an orthonormal basis for its range $\ran{K}\subset H$.
We use its $L\in\Nbb$ leading eigenfunctions
$\{\phi_\ell\}_{\ell=0}^{L-1}$
to define the finite dimensional subspace
$H^L\subset H$ as their linear span.
Next, the density operators and observables introduced
in Section \ref{sec:qmcl-infinite}
are projected from $B(H)$ to $B(H^L)$, which can be
represented by $L\times L$ matrices with respect to our basis.
This process will be presented in Section \ref{sec:qmcl-data}
along with the data based formulation of our closure scheme.

In the next section we describe some of the properties of the
subspace $\ran{K}\subset H$ obtained through the eigendecomposition
of $K$ as they relate to dynamical symmetries.

\subsection{Dynamical symmetries}\label{sec:qmcl-symmetries}
We build the finite dimensional subspace $H^L$ via the
eigenfunctions of the kernel integral operator $K$,
using a kernel function that acts by time delay embedding
each state vector $u$ at a fixed gridpoint $x$.
This choice of kernel function is motivated by the VSA framework
developed in \cite{Giannakis2019vsa}, where the following
material is presented in more detail.

We are particularly interested in PDEs with dynamical symmetries,
which often arise from symmetries on the spatial domain such as
spatial translations or rotations.
We stress that our closure framework does not require the presence
of dynamical symmetries.
However, when symmetries are present, then our choice of kernel
function exhibits some favorable properties.

We consider a symmetry group $G$ acting on the spatial domain $S$
by a continuous left action $\Gamma_g$, which preserves the measure
$\nu$ for all $g\in G$.
In particular, $G$ is the group of dynamical symmetries associated with
the resolved dynamics $\Phi_r^n$, and is generally a subgroup of the
group of symmetries associated with the true dynamics $\Phi^n$.
This is because some of the original symmetries may be lost when
deriving the resolved dynamics by spatial averaging or coarsening.

The actions $\Gamma_g$ induce the associated actions
$\Gamma_{X,g}(u)=u\circ\Gamma_g^{-1}$
on the state space $X$.
Assuming that $X_r$ is invariant under $\Gamma_g$,
they also induce the corresponding actions
$\Gamma_{r,g}(\hat{u})=\hat{u}\circ\Gamma_g^{-1}$
on the resolved state space $X_r$.
The induced actions represent dynamical symmetries of the
resolved dynamics $\Phi_r^n$, which means that they satisfy
the equivariance property
\begin{equation}\label{eq:G-equivariance}
    \Phi_r^n\circ\Gamma_{X,g}=\Gamma_{r,g}\circ\Phi_r^n
\end{equation}
for all $n\in\Nbb$ and $g\in G$.
In addition, we assume that the resolved dynamics is derived
in such a way that the true flux functions $f_j\in H$
satisfy the analogous equivariance property
$f_j(\Gamma_{X,g}(u),x)=f_j(u,\Gamma_g(x))$
for all $x\in S$, $u\in X$, $g\in G$ and $j\in\{0,\ldots,d-1\}$.
As explained below, our choice of kernel naturally exploits the
equivariance property \eqref{eq:G-equivariance}
to bestow an analogous invariance property on the eigenfunctions
of the kernel integral operator.

As an example, the dynamics generated by the advection equation
\eqref{eq:model-pde} with periodic boundary conditions is
equivariant under the continuous group of spatial translations.
The translation symmetry is preserved during the spatial
averaging employed to derive the averaged equation
\eqref{eq:model-pde-coarse}.
As a result, in this case both the true and the resolved
(averaged) dynamics satisfy the same group of dynamical symmetries.
Furthermore, in Section \ref{sec:swe} we will consider the shallow
water equations \eqref{eq:swe} with periodic boundary conditions,
which are also equivariant under spatial translations.
More specifically, we will consider the discretized equations on
a fine grid as the true dynamics, and an averaged version on a
coarse grid as the resolved dynamics.
In that case, the true dynamics will be equivariant under the
discrete group of translations associated with the fine grid.
However, the symmetries of the resolved dynamics will be a
subgroup of that group, containing only the spatial translations
on the coarser grid.

Under the action of map $W_Q$ from \eqref{eq:delay-map},
the product space $\Omega$ is factored into a space of
equivalence classes
$[\omega]_Q\subset\Omega$, $\omega\in\Omega$,
with each class defined as the subset of $\Omega$
on which $W_Q$ is constant.
Namely, a product state $\omega'$ is a member of $[\omega]_Q$
if and only if $W_Q(\omega)=W_Q(\omega')$.
For our choice of map $W_Q$, each equivalence class consists of
the points in $\Omega$ whose resolved dynamical evolution over
$Q$ delays is identical.
From the equivariance property \eqref{eq:G-equivariance}
it follows that, for any $\omega\in\Omega$, the set of all
symmetry transformations of $\omega$ under $G$ is a subset
of its equivalence class $[\omega]_Q$.
This is because any product state vectors that differ only by
a symmetry transformation will demonstrate the same evolution
under the resolved dynamics $\Phi_r^n$, by virtue of
\eqref{eq:G-equivariance}.

Because the kernel function $\kappa$ in \eqref{eq:integral-op}
acts by composition with the delay embedding map $W_Q$,
it follows that $\kappa$ is constant on the $[\cdot]_Q$
equivalence classes by construction.
As a result, every function $f\in\ran{K}$ is constant on the
equivalence classes; namely,
$f(\omega')=f(\omega)$ for all $\omega'\in [\omega]_Q$.
Since each equivalence class $[\omega]_Q$ contains all symmetry
transformations of $\omega$ under $G$, it follows that functions
$f\in\ran{K}$ satisfy the invariance property
\begin{equation}\label{eq:G-invariance}
    f\circ\Gamma_{\Omega,g}=f
\end{equation}
for all $g\in G$, with
$\Gamma_{\Omega,g}=\Gamma_{X,g}\otimes\Gamma_g$
the induced group action on $\Omega$.
In particular, this implies that the eigenfunctions of the integral
operator $K$ satisfy this property as well.
Combined with the fact that all $f\in\ran{K}$ are continuous by
the uniform continuity of $\kappa$,
we conclude that $\ran{K}$ consists entirely of continuous functions
that are constant on the equivalence classes $[\cdot]_Q$.

As a result of \eqref{eq:G-invariance},
by using the eigenfunctions of $K$ we derive a basis of
$\ran{K}$ that is invariant under the actions of $G$,
instead of getting multiple ``copies'' of the same underlying
functions generated by the action of different $g\in G$.
Most importantly for our purposes, this means that fewer eigenfunctions
are needed to represent complex spatiotemporal dynamics,
thereby producing a more efficient compression of the dynamics
from $H$ to subspace $H^L$.

\subsection{Data based formulation}\label{sec:qmcl-data}
In this section we consider the data based formulation
of the extended QMCl scheme developed in the previous sections,
using a single dynamical trajectory to generate our training data.
In the next section we extend this construction to the case of
training data consisting of multiple trajectories.

In the theoretical formulation of QMCl developed above, the state
space $X$ is generally an infinite dimensional Hilbert space of functions
on the spatial domain $S$, with the dynamics map $\Phi^n$ generated
by a PDE.
To be able to simulate the ``true'' dynamics numerically, we
approximate the true dynamics by a discretization of the original PDE
associated with a discretized spatial domain $S_M\subset S$, where
$M\in\Nbb$ denotes the number of gridpoints.
In what follows, we denote by $\nu_M$ a restriction of the Lebesgue
measure $\nu$ to a $\sigma$-algebra of sets generated by the
discretization cells of the spatial domain $S_M$.

With respect to our motivating example of the advection PDE
\eqref{eq:model-pde}, the above corresponds to considering the
true discrete time dynamics $\Phi^n$ as being generated by a
high fidelity discretized version of \eqref{eq:model-pde}.
The resolved dynamics and associated resolved state space are again
defined by forming a spatially averaged version of the
discretized true dynamics.
A detailed description of this construction will be presented in
Section \ref{sec:swe}, where we consider an application of QMCl
to the shallow water equations.

We consider a finite dynamical trajectory of $N\in\Nbb$ samples
$X_N=\{u_n\}_{n=0}^{N-1}\subset X$
with $u_n=\Phi^n(u_0)$, $u_0\in X$,
and the associated sampling measure
$\mu_N=\sum_{n=0}^{N-1}\delta_{u_n}/N$.
The space of observables is defined as the $NM$-dimensional
Hilbert space $H_N=L^2(\Omega,\sigma_{N})$
with product measure $\sigma_N=\mu_N\times\nu_M$.
To construct our closure model we require access to the resolved
state samples $\hat{u}_n=\proj_r(u_n)\in X_r$
and associated flux samples $f_j(u_n,\cdot)\in L^2(S,\nu_m)$,
$j\in\{0,\ldots,d-1\}$.
We denote the product state samples by
$\omega_{nm}=(u_n,x_m)\in\Omega_{NM}$
with $u_n\in X_N$ and $x_m\in S_M$.
We operate under the assumption that the sampling measures
$\mu_N$ converge weakly to the invariant measure $\mu$, in the sense
that $\lim_{N\to\infty}\int_Xf\,d\mu_N=\int_Xf\,d\mu$
for all $f\in C_b(X)$
\cite{Freeman2024}.

In the discrete setting the kernel integral operator
$K_N\colon H_N\to H_N$
corresponding to \eqref{eq:integral-op}
is defined as
\begin{equation}\label{eq:data-integral-op}
    K_Nf=\int_\Omega k(W_Q(\cdot),W_Q(\omega))
        f(\omega)d\sigma_N(\omega)
    =\frac{1}{NM}\sum_{n=0}^{N-1}\sum_{m=0}^{M-1}
        k(W_Q(\cdot),W_Q(\omega_{nm}))f(\omega_{nm}).
\end{equation}
Using the $L$ leading eigenfunctions
$\{\phi_\ell\}_{\ell=0}^{L-1}$ of $K_N$
we construct the $L$-dimensional subspace
$H_N^L\subset H_N$ as their linear span.
Computing the eigenfunctions of $K_N$ with a direct method has
a computational cost that scales with $(NM)^3$, making the computation
prohibitively expensive as $NM$ grows.
To make the computation more tractable, in our numerical experiments
we build a low rank approximation of the kernel matrix $\matr{K}_N$
associated with $K_N$ and use it to compute its eigenfunctions with
reduced cost (Section \ref{sec:swe-numerics}).
In what follows, we make use of the orthogonal projection
$\proj_L\colon H_N\to H_N^L\subset H_N$
to project an operator $A\in B(H_N)$
to $\proj_L A\proj_L\in B(H_N^L)$,
which can then be represented by an $L\times L$
matrix of coefficients with respect to our basis.

Each flux function $f_j\in H_N$ corresponds to a vector in $\Rbb^{NM}$
with entries $f_j(u_n,x_m)$, where $u_n\in X_N$ and $x_m\in S_M$.
Every such function is used to define a multiplication operator
$\tilde{A}_j\in B(H_N)$ and, after orthogonal projection, form
a quantum observable $A_j\in B(H_N^L)$,
which is generally no longer a multiplication operator.
Using our orthonormal basis, each $A_j$ is represented
by an $L\times L$ matrix $\matr{A}_j$ with entries
$A_{j,ik}=\langle\phi_i,\tilde{A}_j\phi_k\rangle_N$.

The discretization process employed by QMCl to generate the
observables $A_j$ is positivity preserving.
This means that if $f_j$ is a sign definite function
(positive or negative), then so will be the resulting observable
$A_j$, producing surrogate fluxes that agree in sign with the
true fluxes.
This is one of the main benefits of mapping classical observables
and density functions to their associated operators before discretizing
them with respect to our basis.
This is in contrast to discretizing the functions directly, which
does not guarantee the preservation of their positivity,
potentially leading to surrogate fluxes that differ in sign
from the true ones.
Of course, this benefit is accompanied by an increase in
computational cost, since every observable is represented by
a matrix of coefficients with respect to our basis, instead
of a vector of coefficients.

To model the unresolved degrees of freedom we form a field of
quantum density operators $\rho\colon S_M\to B(H_N^L)$,
consisting of one density operator $\rho_m=\rho(x_m)\in B(H_N^L)$
for each gridpoint $x_m\in S_M$.
Each density operator $\rho_m$ is constructed as the projection
along a real valued function $h_m\in H_N^L$,
$\rho_m=\langle h_m,\cdot\rangle_Nh_m$.
Taking advantage of its rank-one structure, each $\rho_m$
is represented by the coefficient vector
$\vect{\rho}_m\in\Rbb^L$
with entries $\rho_{m,\ell}=\langle\phi_\ell,h_m\rangle_N$.

Having formed the quantum observables and density operators,
we are now in a position to formulate the calculation of the
required surrogate fluxes.
For each $j\in\{0,\ldots,d-1\}$ we define the surrogate flux
function $\tilde{f}_j\colon S_M\to\Rbb$,
which uses the quantum observable $A_j$
and field of density operators $\rho$
to produce the surrogate fluxes
$\tilde{f}_j(x_m)=\tr(\rho(x_m)A_j)\in\Rbb$
at each gridpoint $x_m\in S_M$.
Using our basis representation this corresponds to
$\tilde{f}_j(x_m)=\vect{\rho}_m^\top\matr{A}_j\vect{\rho}_m$.

The time evolution of the quantum density operators under the dynamics
is carried out using the projected transfer operator
$P_N^L\colon B(H_N^L)\to B(H_N^L)$.
To build $P_N^L$ we begin by defining the shift operator
$P_N\colon H_N\to H_N$
\begin{equation*}
    P_Nf(\hat{u}_n,x_m)=\begin{cases}
        0& \text{if $n=0$}\\
        f(\hat{u}_{n-1},x_m)& \text{if $1\leq n<N$}.
        \end{cases}
\end{equation*}
After orthogonal projection we obtain $P_N^L$
and its $L\times L$ coefficients matrix $\matr{P}_N^L$.
Given a quantum density operator $\rho_n(x_m)$
at time $t_n$ with vector representation $\vect{\rho}_n$
we obtain the updated state $\vect{\rho}_{n+1}$
by
\begin{equation*}
    \vect{\rho}_{n+1}=\frac{\matr{P}_N^L\,\vect{\rho}_n}
        {\lVert\matr{P}_N^L\,\vect{\rho}_n\rVert_2}
\end{equation*}
where the renormalization is performed because the matrix
$\matr{P}_N^L$ is generally not orthogonal.

Finally, we arrive at the bayesian conditioning of quantum density
operators.
Each density operator $\rho_m$ is conditioned by a quantum effect
operator $e_m\in B(H_N^L)$ formed based on a given product state
vector $\omega_m=(u,x_m)$.
To build $e_m$ from $\omega_m$ we use the conditioning kernel function
$\kappa_c\colon\Omega\times\Omega\to\Rbb$,
$\kappa_c(\omega,\omega')=k_c(W_J(\omega),W_J(\omega'))$
which acts by composition of a gaussian kernel function
$k_c\colon\Rbb^J\times\Rbb^J\to\Rbb$
and an embedding map
$W_J\colon\Omega\to\Rbb^J$
with parameter $J\in\Nbb$.
The map $W_J$ returns the resolved state component of the
$J$ nearest spatial neighbors of $\omega_m$ in a specified order.
For example, for a one dimensional grid and odd $J$ its action is
given by
\begin{equation*}
    W_J((u,x_m))=\bigl(\proj_r(u)(x_{m-(J-1)/2}),\,
        \ldots,\,\proj_r(u)(x_m),\,\ldots,\,
        \proj_r(u)(x_{m+(J-1)/2})\bigr).
\end{equation*}
This can be generalized analogously for higher dimensional grids.
The gaussian kernel $k_c$ then acts by
\begin{equation}\label{eq:data-kernel-cond}
    k_c(W_J(\omega),W_J(\omega'))=
        \exp\bigl(-\frac{1}{\epsilon' J}
        \lVert W_J(\omega)-W_J(\omega')\rVert^2\bigr)
\end{equation}
where $\lVert\cdot\rVert$ denotes the euclidean norm in $\Rbb^J$
and $\epsilon'>0$ a bandwidth parameter.
Note that the bandwidth parameter $\epsilon'$ is different
from the one introduced earlier in kernel
\eqref{eq:data-kernel};
it can again be selected using the tuning procedure developed in
\cite{Coifman2008}.
In our numerical experiments we employ the above kernel using a
variable bandwidth; details on this construction and an algorithm
for computing the variable bandwidth can be found in
\cite{Berry2016,Giannakis2019acha}.

For the given product state vector $\omega_m$, we use the
kernel $\kappa_c$ and the full training dataset $\Omega_{NM}$
to produce the function $f_m\in H_N$,
$f_m(\cdot)=\kappa_c(\omega_m,\cdot)$,
which compares the vector $\omega_m$
with every other product state vector included in our training dataset.
Intuitively, the function $f_m$ can be thought of as a similarity function
on the dataset $\Omega_{NM}$, measuring how ``similar'' the given state
$\omega_m$ is to the states included in $\Omega_{NM}$.
Because kernel $\kappa_c$ acts by composition with $W_J$,
it compares the nearest spatial neighbors of $\omega_m$
with the nearest neighbors of every other available state vector,
mimicking the structure of the stencil used for the spatial
discretization of the considered dynamics.
This makes for a stronger similarity condition than comparing values
of states at only a single gridpoint, forcing the kernel to be more
selective in the states it considers to be dynamically similar to
the given $\omega_m$.

Next, we use the square root $f_m^{1/2}$ to form a multiplication
operator and, after orthogonal projection, the desired quantum
effect $e_m\in B(H_N^L)$.
This is then used to condition the prior density operator
$\rho_m$ via the quantum Bayes rule \eqref{eq:quantum-bayes},
thereby creating the posterior density operator
$\rho_m\vert_{e_m}$ that encodes the state space similarities
captured by $f_m$.
Based on our basis representation, this operation takes the form
\begin{equation*}
    \vect{\rho}_m\vert_{e_m}
        =\frac{\matr{e}_m\vect{\rho}_m}
        {\lVert\matr{e}_m\vect{\rho}_m\rVert_2}
\end{equation*}
where $\matr{e}_m$ denotes the $L\times L$ matrix of coefficients
of effect $e_m$.
We refer to the map $\omega_m\mapsto e_m$
as the \emph{feature map} and to $f_m$ as the
\emph{feature vector} produced by kernel $\kappa_c$
from the product state vector $\omega_m$.
Explanation for why we use the square root of the feature vector
$f_m^{1/2}$ instead of $f_m$ to build $e_m$ is given in
\cite{Freeman2024}.

Using the procedure outlined above, the density operator
$\rho_m=\rho(x_m)$ associated with every point $x_m$ on
our grid can be conditioned based on the corresponding
product state vector $(u,x_m)$.
That is, given a new state vector $u\in X$, the density operator
$\rho_m$ is conditioned based on the spatially local information
$\omega_m=(u,x_m)$, by comparing this product state vector to every
other product state vector available in our training dataset.
Because the spatially local value $\omega_m$ is used for the
conditioning process (rather than the full state $u$),
spatial translations of the state $u$ are naturally factored during
the conditioning process.
This means that if a spatial translation of the state $u$ is present
in our training dataset, then the similarity between the product state
$\omega_m$ and the corresponding value of the spatial translation of
$u$ will be identified in the conditioning process.
For PDEs with spatial translation symmetry, this means that we do not
have to artificially inflate our training dataset by including all
spatial translations of the sampled states.

\subsection{Multi-trajectory data}\label{sec:qmcl-multi}
In the previous section we presented the data based implementation
of our closure scheme using a single dynamical trajectory as our
training dataset.
In this section we generalize this construction to the case where
the training dataset consists of multiple trajectories.
We begin by motivating the use of multiple trajectories,
before presenting the modifications required in the
theoretical and data based formulations of the scheme.

Using a single dynamical trajectory is generally sufficient when
considering a dynamical system possessing a unique physical measure,
namely a measure that can be sampled by a set of
initial conditions of positive measure with respect to the
ambient measure \cite{LSYoung2002,Blank2017}.
Using multiple trajectories is of interest when considering
dynamical systems possessing multiple physical measures;
for instance, hamiltonian systems or hyperbolic systems with
several competing basins of attraction.
Additionally, using multiple trajectories is necessary in the
forecasting of parameter dependent systems, whose qualitative
dynamics changes with respect to a bifurcation parameter.
Finally, even in systems with a unique physical measure,
being able to use multiple trajectories allows one to consider
optimal sampling strategies.
In designing such a strategy one tries to optimize the sampling
rate with respect to a chosen objective, instead of relying on the
inherent rate provided by a single trajectory.

For example, in Section \ref{sec:swe} we will consider an
application of our scheme to a closure problem for the shallow
water equations \eqref{eq:swe}, a hyperbolic, energy conservative
system of PDEs.
Depending on the initial condition, the solution trajectory of
equations \eqref{eq:swe} can reach different limit sets, since there
is no dissipation to focus the global dynamics to a lower dimensional
subset of its state space.
As such, using several trajectories becomes important, allowing
us to sample different limit sets of the dynamics.

We consider a finite collection of $I\in\Nbb$ invariant probability
measures $\{\mu_i\}$ on $X$ with mutually disjoint supports,
$i\in\{0,\ldots,I-1\}$,
and define the invariant probability joint measure
$\mu=\sum_{i=0}^{I-1}\mu_i/I$
as their average.
By construction, the Hilbert space of observables
$H=L^2(\Omega,\mu\times\nu_M)$
admits the decomposition $H\simeq\oplus_{i=0}^{I-1}H_i$
with component spaces $H_i=L^2(\Omega,\mu_i\times\nu_M)$,
where the isomorphism is up to a rescaling of the inner products
of $H_i$ by the factor $1/I$.
Our data based formulation mirrors the one presented in
Section \ref{sec:qmcl-data},
while also taking advantage of the decomposition of $H$ into its
component subspaces.

We generate the finite sampled trajectories
$\{u_n^{(i)}\}_{n=0}^{N_i-1}$ of $N_i\in\Nbb$ samples each,
where $i\in\{0,\ldots,I-1\}$ as above.
Each trajectory gives rise to its associated empirical sampling measure
$\mu_{i,N_i}=\sum_{n=0}^{N_i-1}\delta_{u_n^{(i)}}/N_i$,
where we assume that each sampling measure $\mu_{i,N_i}$
converges to $\mu_i$ in the weak-$*$ sense as $N_i\to\infty$.
Denoting $N=\sum_{i=0}^{I-1}N_i$,
we define the Hilbert space
$H_N=L^2(\Omega,\mu_N\times\nu_M)$
with the joint sampling measure
$\mu_N=\sum_{i=0}^{I-1}\mu_{i,N_i}/I$.
Since the individual trajectories are assumed to be mutually
disjoint, we have the analogous decomposition
$H_N\simeq\oplus_{i=0}^{I-1}H_{N_i}$.

Leveraging the decomposition of $H_N$, we use each individual
trajectory and the approach described in
Section \ref{sec:qmcl-data}
to build an eigenfunction basis for each component space
$H_{N_i}$ separately.
Relying on the $L_i\in\Nbb$ leading eigenfunctions, we then form the
subspace $H_{N_i}^{L_i}\subset H_{N_i}$.
Denoting $L=\sum_{i=0}^{I-1}L_i$,
we define the combined space
$H_N^L=\oplus_{i=0}^{I-1}H_{N_i}^{L_i}$,
with the direct sum of the component bases producing a basis
for the combined space.

The basis of $H_N^L$ derived in this way consists of a total of $L$
elements, each being an eigenfunction of a block diagonal operator
$K_N=K_{N_0}\oplus\ldots\oplus K_{N_{I-1}}$
with the component kernel integral operators
$K_{N_i}$ in its diagonal blocks.
More specifically, each component basis vector $\phi_\ell^{(i)}$
represents the values of a function on $N$ samples,
with all of its values being zero except for the block of values
that corresponds to the $N_i$ samples of the $i$th trajectory it
was derived from.
This approach has lower computational cost, as it requires that we
solve a total of $I$ eigenvalue problems for an $N_i\times N_i$
matrix each, instead of a single eigenvalue problem for an
$N\times N$ matrix.
However, our approach also requires that each component measure
$\mu_i$ is individually adequately sampled,
as opposed to adequate joint sampling of measure $\mu$.

Having derived the basis for space $H_N^L$,
the construction of the rest of the scheme proceeds as in
Section \ref{sec:qmcl-data}.
The only additional difference is in the definition of the
shift operator $P_N$,
where time shifts are now applied only \emph{within} each
trajectory, and the first temporal sample of each trajectory
is mapped to zero.
Namely, we make sure not to shift samples
\emph{across} different trajectories.

After forming our QMCl closure model in $H_N^L$,
our goal is to use it to predict the dynamics of a chosen
dynamical system for out of sample initial conditions;
namely, initial conditions lying outside of the basins of the
component measures $\mu_i$ sampled in the training stage.

\section{Application}\label{sec:swe}
The shallow water equations govern the motion of an incompressible
fluid when its vertical length scale can be considered negligible
compared to its horizontal length scale.
The equations are derived from the conservation laws for mass and
momentum, with various forms finding applications in areas such as
ocean dynamics, often incorporating additional terms modeling
the effects of bottom topography and Coriolis force
\cite{Vallis2017,Whitham1999}.

We consider the shallow water equations on a bounded one dimensional
spatial domain $S$ with periodic boundary conditions.
The system of equations governs the evolution of the water height
$h$ and momentum $q=hv$ fields, with $v$ denoting the velocity field.
The governing equations read
\begin{equation}\label{eq:swe}
\begin{aligned}
    \partial_th+\partial_xq&=0\\
    \partial_tq+\partial_x(q^2/h+1/2Fr^{-2}h^2)&=0
\end{aligned}
\end{equation}
subject to an appropriate initial condition.
In the above $\partial_t$ and $\partial_x$ denote differentiation
with respect to time and space respectively, while $Fr$ denotes
the flow's Froude number, which is a dimensionless ratio of inertial
over gravitational forces.
Using the state variable
$u=(h,q)\in X\subset L^2(S,\nu;\Rbb^2)$
and function
$f(u)=(q,q^2/h+1/2Fr^{-2}h^2)$
we rewrite \eqref{eq:swe} in the compact form
\begin{equation}\label{eq:swe2}
    \partial_tu+\partial_xf(u)=0.
\end{equation}

In addition to its applications in ocean dynamics,
the shallow water equations problem \eqref{eq:swe2}
has additional properties that make it a good test case
for our QMCl closure framework.
First, the dynamics governed by \eqref{eq:swe2} is equivariant
under the actions of the spatial translation group on $S$.
Second, forming a training dataset that is sufficiently representative
of a range of its dynamics requires the use of multiple dynamical
trajectories.
Third, previous works have considered various closure schemes for
variations of \eqref{eq:swe2}, demonstrating the amenability
of the equations to closure problems that arise from spatial
coarsening \cite{Zacharuk2018,Timofeyev2025}.
In particular, our discretization and coarsening methodology
is based on \cite{Timofeyev2025}.
Finally, the equations act as a good idealized model for more
complex nonlinear hyperbolic PDEs that are prevalent in fluid
and plasma dynamics modeling.

\subsection{Discretization}\label{sec:swe-discretization}
We employ a finite volume spatial discretization of
\eqref{eq:swe2}
using a uniform fine grid of $M_f\in\Nbb$ cells with length
$\Delta x=L/M_f$ each, centered at points
$x_m\in S$, $m\in\{0,\ldots,M_f-1\}$.
We use $S_f=\{x_m\}_{m=0}^{M_f-1}\subset M$
to denote the center points of the fine discretization cells.
The cell averaged state variables $u_m$, $m\in\{0,\ldots,M_f-1\}$,
are defined by taking the spatial average of the original
state function $u$ over the cell centered at $x_m$.
Without misunderstanding, we reuse the same symbol $u$ to
denote the averaged state variables, which can be represented
by the vector $\vect{u}\in\Rbb^{2M_f}$ for every time $t\geq 0$.

The temporal evolution of the fine state variables $u_m$
is governed by the system of equations
\begin{equation}\label{eq:swe-semi}
    \dot{u}_m=\frac{1}{\Delta x}(F_m-F_{m+1})
\end{equation}
using the local Lax-Friedrichs (LLF) numerical flux
\begin{equation}\label{eq:llf-flux}
    F_m=\frac{1}{2}[f(u_{m-1})+f(u_m)]
        -\frac{\lambda_m}{2}(u_m-u_{m-1})
\end{equation}
with
\begin{equation}\label{eq:llf-wavespeed}
    \lambda_m=\max\,\bigl(\lvert q_{m-1}/h_{m-1}\rvert+Fr^{-1}\sqrt{h_{m-1}},
        \,\vert q_m/h_m\rvert+Fr^{-1}\sqrt{h_m}\bigr).
\end{equation}
In the above, subscripts denote components of the associated vectors
while $F_m$ denotes the numerical flux at the left face of the cell
centered at $x_m\in S_f$ \cite{LeVeque2002}.
Based on \eqref{eq:llf-flux} and \eqref{eq:llf-wavespeed}
we introduce the notation
\begin{equation*}
    F_m=\Fcal_{LLF}(u_{m-1},u_m)
\end{equation*}
which will be used later on to define the coarse grid fluxes.

We now form a coarse grid and associated coarse state variables
defined as local spatial averages of the fine variables $u_m$.
In particular we form $M_c\in\Nbb$ coarse grid cells by merging every
$M_s\in\Nbb$ fine cells; in this way, $M_c=M_f/M_s$ and each coarse
cell has length $\Delta\hat{x}=M_s\Delta x$.
We define the resulting coarse grid
$S_c=\{\hat{x}_m\}_{m=0}^{M_c-1}\subset S$,
$m\in\{0,\ldots,M_c-1\}$,
representing the center points of the coarse grid cells.
Each coarse grid state variable $\hat{u}_m$ is defined as
the average of the fine state variables included within its
associated coarse grid cell.
The coarse variables are represented by the vector
$\vect{\hat{u}}\in\Rbb^{2M_c}$
for every time $t\geq 0$.
Their temporal evolution is governed by the system of equations
\begin{equation}\label{eq:swe-semi-coarse}
    \dot{\hat{u}}_m=\frac{1}{\Delta\hat{x}}
        [(\hat{F}_m+G_m)-(\hat{F}_{m+1}+G_{m+1})]
\end{equation}
with coarse grid fluxes
\begin{equation}\label{eq:coarse-flux}
    \hat{F}_m=\mathcal{F}_{LLF}(\hat{u}_{m-1},\hat{u}_m)
\end{equation}
and subgrid fluxes
\begin{equation}\label{eq:subgrid-flux}
    G_m=F_{M_sm}-\hat{F}_m.
\end{equation}
The system \eqref{eq:swe-semi-coarse} is equivariant under the
actions of the discrete group of spatial translations on $S_c$.

By definition, the coarse grid fluxes \eqref{eq:coarse-flux}
depend only on the coarse variables.
On the contrary, the subgrid fluxes \eqref{eq:subgrid-flux}
depend on both the coarse and the original fine variables;
in particular, they depend on the values of $u_m$
in the fine cells that border the interfaces of each coarse cell.
For the physical system considered, the subgrid fluxes act as
antidiffusive corrections to the coarse grid fluxes, increasing the
accuracy of the discretization compared to one performed entirely
on the coarse grid \cite{Timofeyev2025}.

The temporal discretization of evolution equations
\eqref{eq:swe-semi} and \eqref{eq:swe-semi-coarse}
is carried out using the modified Euler method
(two-stage Runge-Kutta), which has favorable structure
preservation properties for the present problem
\cite{Timofeyev2025}.
In what follows, we consider the coarse grid variables as the
resolved degrees of freedom and the discrete time version of
\eqref{eq:swe-semi-coarse}
as the equation generating the resolved dynamics.
To close the system of equations we need to provide a closure model
that can infer the subgrid fluxes $G_m$ from the resolved variables
$\hat{u}_m$, thereby removing the explicit dependence on the
original fine grid variables $u_m$.

\subsection{Closure model}\label{sec:swe-closure}
To generate our training dataset and be able to compare our
surrogate model with a reference model acting as the ground truth,
we approximate the ``true'' shallow water dynamics by the discrete
time version of the fine grid dynamics
\eqref{eq:swe-semi}.
As such, the full state space of the dynamics is a Hilbert subspace
of the set of $L^2$ functions on the grid produced by our fine
discretization.
More specifically we set $X\subset L^2(S,\nu_f;\Rbb^2)$,
where $\nu_f$ denotes the restriction of the Lebesgue volume
measure $\nu$ to the $\sigma$-algebra generated by the fine
discretization cells.

As mentioned in the previous section, we consider the coarse
grid variables as our resolved degrees of freedom.
Therefore, we define the resolved state space
$X_r=L^2(S,\nu_c;\Rbb^2)\cap X$
and decompose the full state space as $X=X_r\oplus X_u$,
with unresolved state space $X_u=X_r^{\perp}$.
In the definition of $X_r$, the measure $\nu_c$ corresponds to
the restriction of $\nu$ to the $\sigma$-algebra of sets generated
by the coarse discretization cells.
As a result, the resolved state space can be considered as a
space of piecewise constant functions, each taking a constant
value on each cell formed by the coarse discretization.
With these choices in place, the full dynamics
$\Phi^{n\Delta t}$
is given by the discrete time version of
\eqref{eq:swe-semi}
and the resolved dynamics
$\Phi_r^{n\Delta\hat{t}}$
by the discrete time version of
\eqref{eq:swe-semi-coarse},
with $n\in\Nbb$ and respective timesteps
$\Delta t>0$ and $\Delta\hat{t}=M_s\Delta t$.

We define the product state space $\Omega=X\times S$
and the Hilbert space of observables $H_N=L^2(\Omega,\sigma_N)$,
with total number of samples $N\in\Nbb$
and measures $\sigma_{N}=\mu_N\times\nu_M$
and $\nu_M=\nu_c$.
For notational consistency with the previous sections,
from now on we denote by $M=M_c$ the number of coarse grid cells
and by $S_M=S_c\subset S$ the set of gridpoints representing the
coarse grid cells.
We introduce two flux functions
$f_j\in H_N$, $j\in\{0,1\}$,
producing the fluxes for the height and momentum
variables respectively,
and denote their collection by $f=\{f_j\}_{j=0}^1$.
For every state $u\in X$ the associated fluxes $f_j(u)$
are given by \eqref{eq:subgrid-flux}.

Using the full dynamics \eqref{eq:swe-semi}
we generate $N$ samples of the resolved state variables
$\{\hat{u}_n\}_{n=0}^{N-1}$
and the fluxes $\{f(u_n,\cdot)\}_{n=0}^{N-1}$,
which will be used as our training data and may consist
of one or more dynamical trajectories.
Each sample $\hat{u}_n$ or $f(u_n,\cdot)$
is represented by a vector in $\Rbb^{2M}$
holding the values corresponding to the height $h$
and momentum $q$ variables in each coarse grid cell.

To form our closure model we define a field of quantum
density operators $\rho\colon S_M\to B(H_N)$
and two quantum observables $A_j\in B(H_N)$, $j\in\{0,1\}$.
Each density operator $\rho(x_m)$ will be used to model in a
statistical sense the unresolved variables in the coarse grid
cell centered at $x_m\in S_M$.
From the two quantum observables, observable $A_0$ will be used
to produce the surrogate fluxes for the height $h$ variables and
is defined as multiplication by the true flux function $f_0$;
analogously, observable $A_1$ will produce the surrogate fluxes
for the momentum $q$ variables and is defined as multiplication
by $f_1$.
With these definitions in place, the surrogate fluxes are
computed using the two surrogate flux functions
$\tilde{f}_j\colon S_M\to\Rbb$
with action $\tilde{f}_0(x_m)=\tr(\rho(x_m)A_0)$
and $\tilde{f}_1(x_m)=\tr(\rho(x_m)A_1)$
for every $x_m\in S_M$.

In \ref{app:implementation} we list the sequence of steps used
to build the data based QMCl model and apply it to predict the
dynamics based on a given initial condition.

\subsection{Numerical results}\label{sec:swe-numerics}
We consider the shallow water equations \eqref{eq:swe2}
on the spatial domain $S=[-25,25]$
with Froude number $Fr=1/\sqrt{2g}$,
where $g$ denotes the acceleration of gravity.
We perform the spatial discretization using a fine grid of
$M_f=1920$ cells and the temporal discretization using the
midpoint Euler method with timestep $\Delta t=0.1\Delta x$,
where $\Delta x$ is the length of each fine cell.

\subsubsection*{Training}
We define the one-parameter family of initial conditions
\begin{equation}\label{eq:ic-family}
\begin{aligned}
    h_0(x) &= 1 + 0.3(1-\delta/2)\sin\Bigl(\frac{2\pi}{L_S}
        3.5(1-\delta/2)x+\pi/6\Bigr)\\
    v_0(x) &= 1 + 0.2(1-\delta)\sin\Bigl(\frac{2\pi}{L_S}
        3(1-\delta)x\Bigr)
\end{aligned}
\end{equation}
with $x\in S$, $L_S=50$, parameter $\delta\in [0,1]$,
$h_0$ the initial condition for height and
$q_0=h_0v_0$ the initial condition for momentum.
The solutions generated by this family of initial conditions
feature a number of wave fronts traveling in different directions
and interacting with one another.
Figures \ref{fig:res01a} and \ref{fig:res02a}
present space-time heatmaps of the solutions that will be used
to test the predictive performance of our closure model.

\begin{figure}[t]
    \centering
    \includegraphics[width=.7\textwidth]{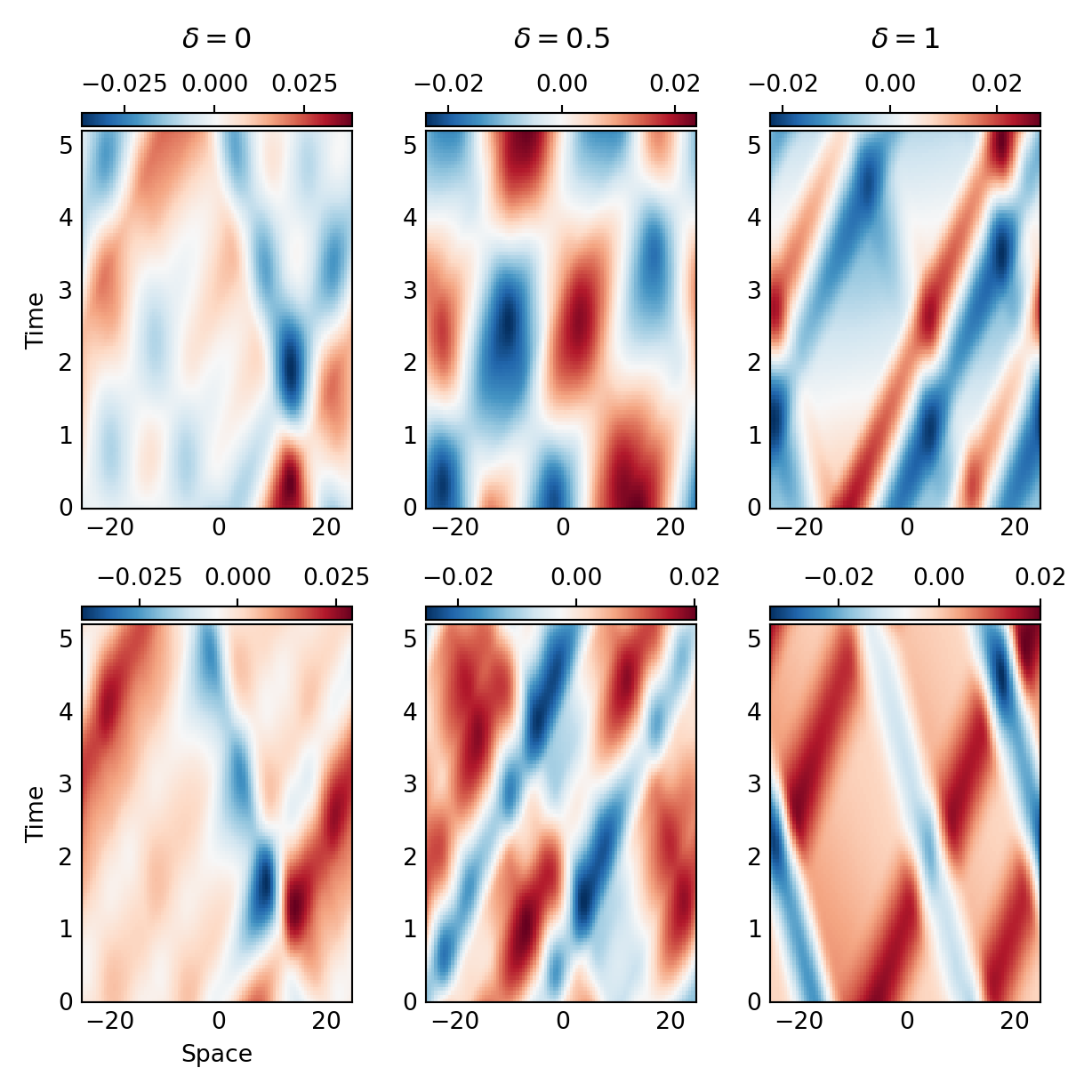}
    \caption{Space-time heatmaps depicting two eigenfunctions
    $\phi_\ell$ of the kernel integral operator $\matr{K}_N$
    from \eqref{eq:data-integral-op}
    for each of the three trajectories used as our training data.
    The top row shows the third eigenfunction for each of the
    three trajectories; the bottom row the fourth eigenfunction
    for each trajectory.
    From left to right, the first column corresponds to the
    trajectory data obtained for $\delta=0$,
    the second to $\delta=0.5$
    and the third to $\delta=1$.}
    \label{fig:vsa-basis}
\end{figure}

Using the assumptions and methodology introduced
in Section \ref{sec:qmcl-multi} for dealing
with training data consisting of multiple trajectories,
we consider $\mu$ to be the invariant measure generated by
the full dynamics \eqref{eq:swe-semi}
restricted to the family of initial conditions
\eqref{eq:ic-family} for $\delta\in [0,1]$.
To sample this measure we use three trajectories generated
using the values $\delta\in\{0,0.5,1\}$.
For each trajectory we integrate the full dynamics
\eqref{eq:swe-semi} in time for 12,200 timesteps and then
sample the next 3,260 timesteps.
Each sample contains the values of functions $h$ and $q$
on the fine grid.

To derive the resolved state samples that will be used as our
training data, we spatially average the collected samples.
By averaging the state data over every $M_s=20$ fine cells
we arrive at a coarse grid of $M=96$ cells of length
$\Delta\hat{x}=20\Delta x$ each.
Given that each cell is now $20$ times longer, we also increase
the sampling timestep by the same amount to set
$\Delta\hat{t}=20\Delta t$.
From the original 3,260 samples per trajectory generated on the
fine grid, we keep every 20th sample, leading to 163 samples
per trajectory.
For each one of those samples we spatially average the solution
over every $20$ cells to derive our resolved state samples,
namely the averaged values of $h$ and $q$ over the coarse grid.
We then derive the subgrid flux training samples using
\eqref{eq:subgrid-flux}.

Using $Q=64$ delays, these samples yield 100 training samples
per trajectory in delay embedded form, adding up to $N=300$
samples in total.
Each resolved state sample consists of $2M=192$ values,
96 for each of the height and momentum variables.
With this training dataset, the dimension of the product
Hilbert space $H_N$ is equal to $NM=28800$.
Employing the approach described in
Section \ref{sec:qmcl-multi}
we build an orthonormal basis $\{\phi_\ell\}_{\ell=0}^{L-1}$
with a total of $L=6144$ eigenfunctions and define
$H_N^L$ as their linear span.
Each eigenfunction is represented by a vector in $\Rbb^{NM}$.
Finally, we set $J=5$ to define the stencil map $W_J$
that will be used for the bayesian conditioning of the quantum
density operators.
The data based realization of the closure model developed
in Section \ref{sec:swe-closure}
will take place in the constructed subspace $H_N^L$.

The number of delays $Q$ was chosen by a combination of
physical reasoning and trial and error to balance the
effectiveness of the derived basis in representing the true
flux samples and the cost associated with its computation.
A useful rule of thumb is to select a delay window that is
comparable to the dominant timescale of the physical processes
of interest.
In our case, we see that the traveling wave fronts interact
approximately every 30 timesteps ($\approx 1.5$ time units)
for the solution in Figure \ref{fig:res01a},
and approximately every 50 timesteps ($\approx 2.5$ time units)
for the solution in Figure \ref{fig:res02a},
with similar results for additional solutions not shown here.
Using $Q=64$ timesteps is sufficient to ensure that at least one
period of this interaction is captured in each delay embedded sample.
Additionally, increasing the value of $Q$ further did not lead
to noticeable differences in our results.
The stencil size $J$ used for the conditioning of the quantum density
operators was chosen empirically.
Increasing the value of $J$ further was shown to slightly improve
the performance of the bayesian conditioning, but with an associated
increase in computational cost.

To form the orthonormal basis
$\{\phi_\ell\}_{\ell=0}^{L-1}$
we must compute the eigenvectors of the $NM\times NM$
symmetric kernel matrix $\matr{K}_N$ representing the operator
\eqref{eq:data-integral-op}.
To reduce the computational cost of this operation
we perform a partial Cholesky factorization of
$\matr{K}_N$ to approximate it by a lower rank matrix
$\tilde{\matr{K}}_N=\matr{F}\matr{F}^T$
with $\matr{F}\in\Rbb^{NM\times r}$ and $r<NM$.
We perform the factorization using the rank parameter
$r=L=6144$ and the randomized algorithm developed in
\cite{YChen2024}.
We then use the approximate kernel matrix to perform its
bistochastic normalization and solve the eigenvalue
problem with reduced cost
\cite{Vales2025evd}.

A selection of the resulting eigenfunctions is presented in
Figure \ref{fig:vsa-basis},
where we show two eigenfunctions for each of the three
trajectories used to build our training dataset.
The shown eigenfunctions exhibit propagating fronts
which resemble the main features of the true state and flux
samples presented in Figures
\ref{fig:res01a} and \ref{fig:res02a}.
More specifically, they are approximately constant on the
characteristic lines of the propagating wave fronts,
as expected from their symmetry invariance property.
Similar features are present in the rest of the eigenfunctions
not shown here.
As explained in Sections
\ref{sec:qmcl-finite} and \ref{sec:qmcl-symmetries},
the strong dynamical relevance of each eigenfunction is a result
of its dynamical symmetry invariance and of the fact that generally
it cannot be represented as the tensor product of a single pair of
spatial and temporal modes.
These properties allow us to compute a compressed representation
of the true dynamics with a moderate number of eigenfunctions,
enhancing the computational efficiency of the resulting closure
scheme.

\subsubsection*{Prediction}
To test the predictive performance of the surrogate model
derived by QMCl we use two initial conditions obtained from
\eqref{eq:ic-family} for $\delta\in\{0.25,0.75\}$,
which are trajectories not included in the training dataset.
After integrating the true dynamics \eqref{eq:swe-semi}
for 12,200 timesteps we spatially average the found solutions
to generate the initial conditions that will be used to initialize
the resolved state in QMCl.
Using those initial conditions we integrate in time the coarse
dynamics \eqref{eq:swe-semi-coarse} using timestep $\Delta\hat{t}$
and the subgrid fluxes predicted by QMCl.
In what follows, we present results obtained after integrating for
a total of 120 timesteps, while performing the bayesian conditioning
of the quantum density operators every 10 timesteps.

\begin{figure}[t!]
    \centering
    \includegraphics[width=.7\textwidth]{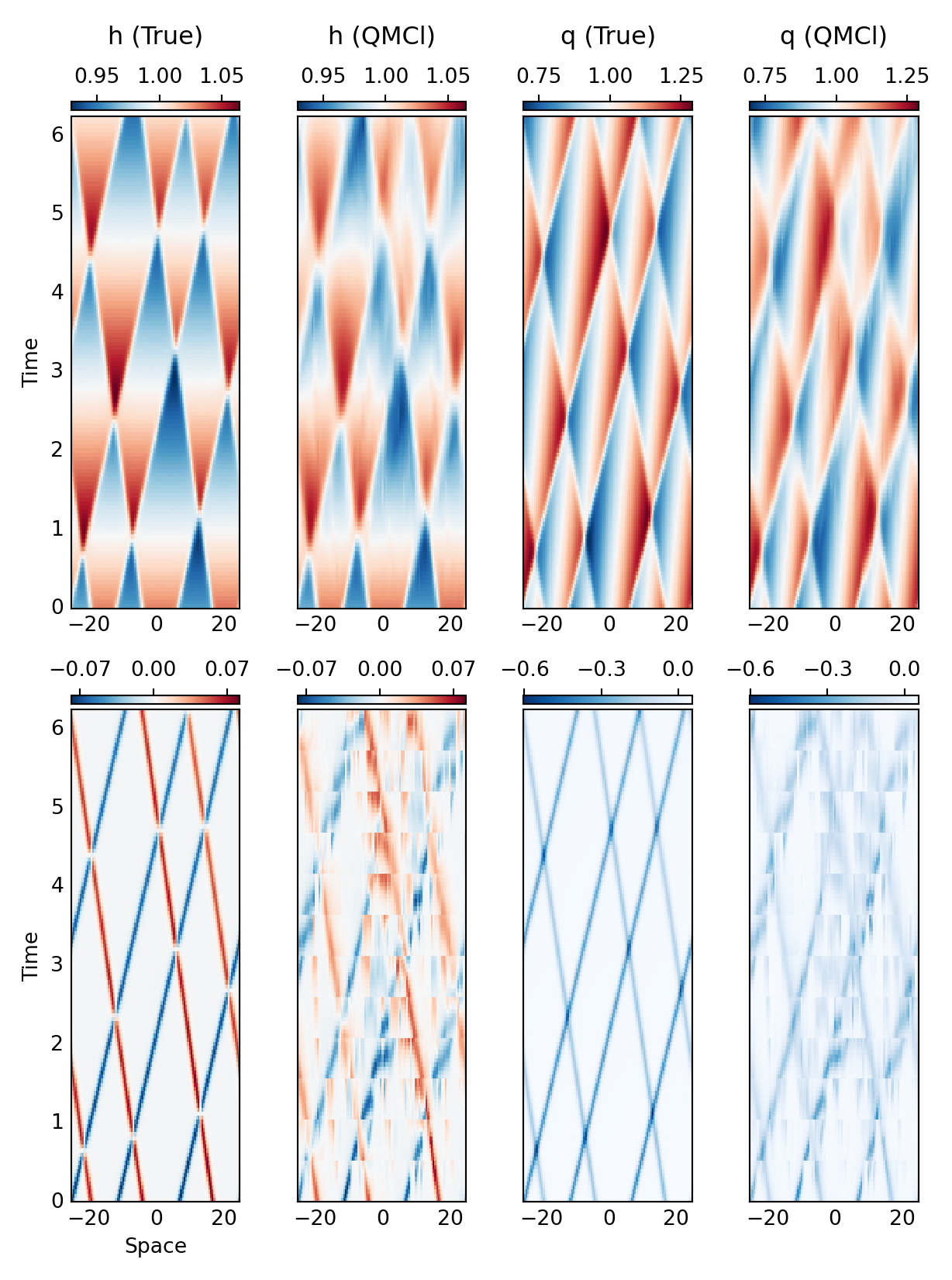}
    \caption{Space-time heatmaps for the trajectory with $\delta=0.25$.
    The top row shows the results for the resolved state variables;
    the bottom row for the subgrid fluxes.
    For each row, from left to right: the first two plots show the results
    for the height field $h$, with the true results on the left
    and the predicted results on the right;
    the next two plots show the results for the momentum field $q$,
    with true results on the left and predicted on the right.}
    \label{fig:res01a}
\end{figure}

\begin{figure}[t!]
    \centering
    \includegraphics[width=.7\textwidth]{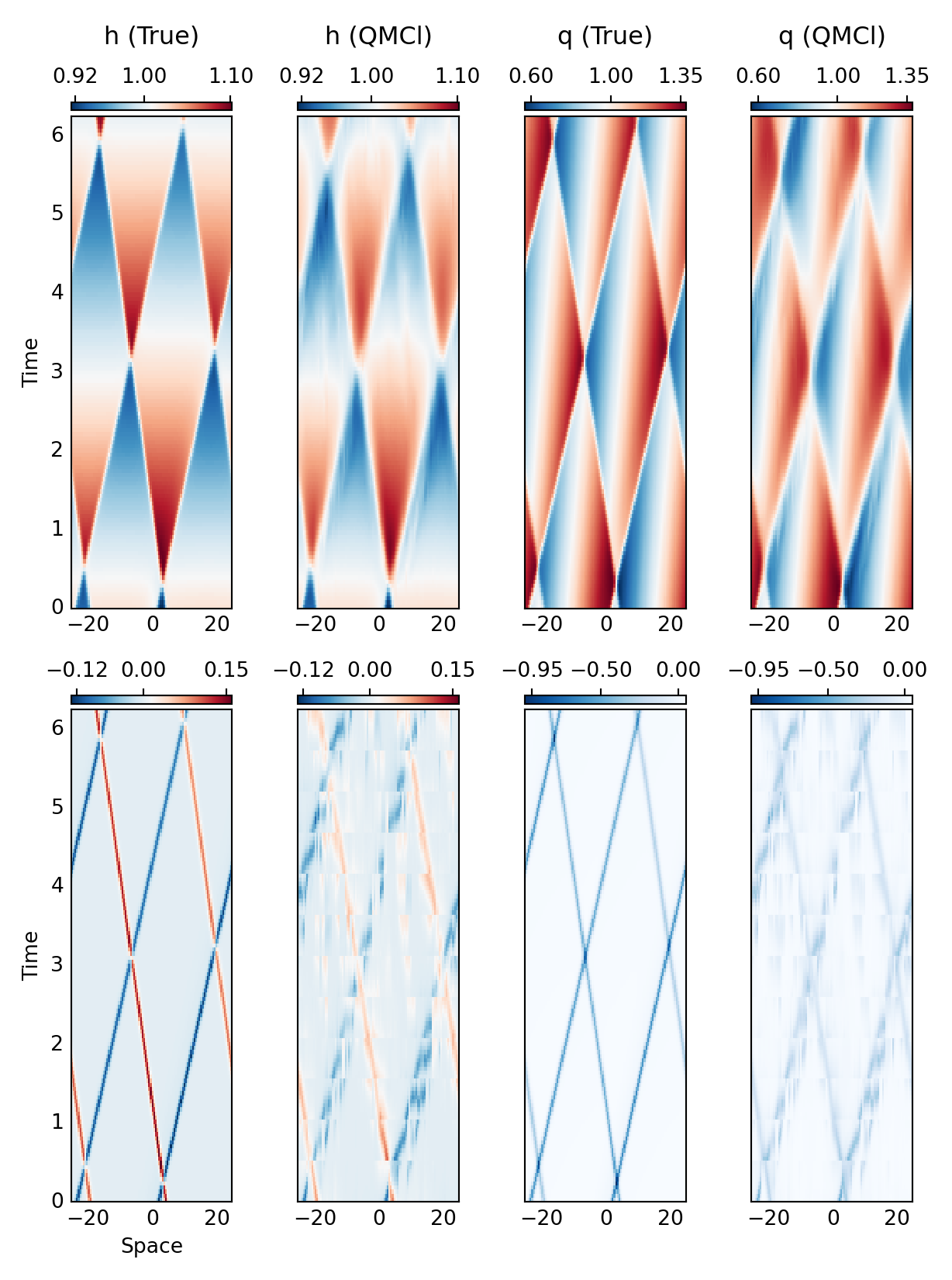}
    \caption{Space-time heatmaps for the trajectory with $\delta=0.75$.
    The layout of the plots is as in Figure \ref{fig:res01a}.}
    \label{fig:res02a}
\end{figure}

\begin{figure}[t]
    \centering
    \includegraphics[width=.7\textwidth]{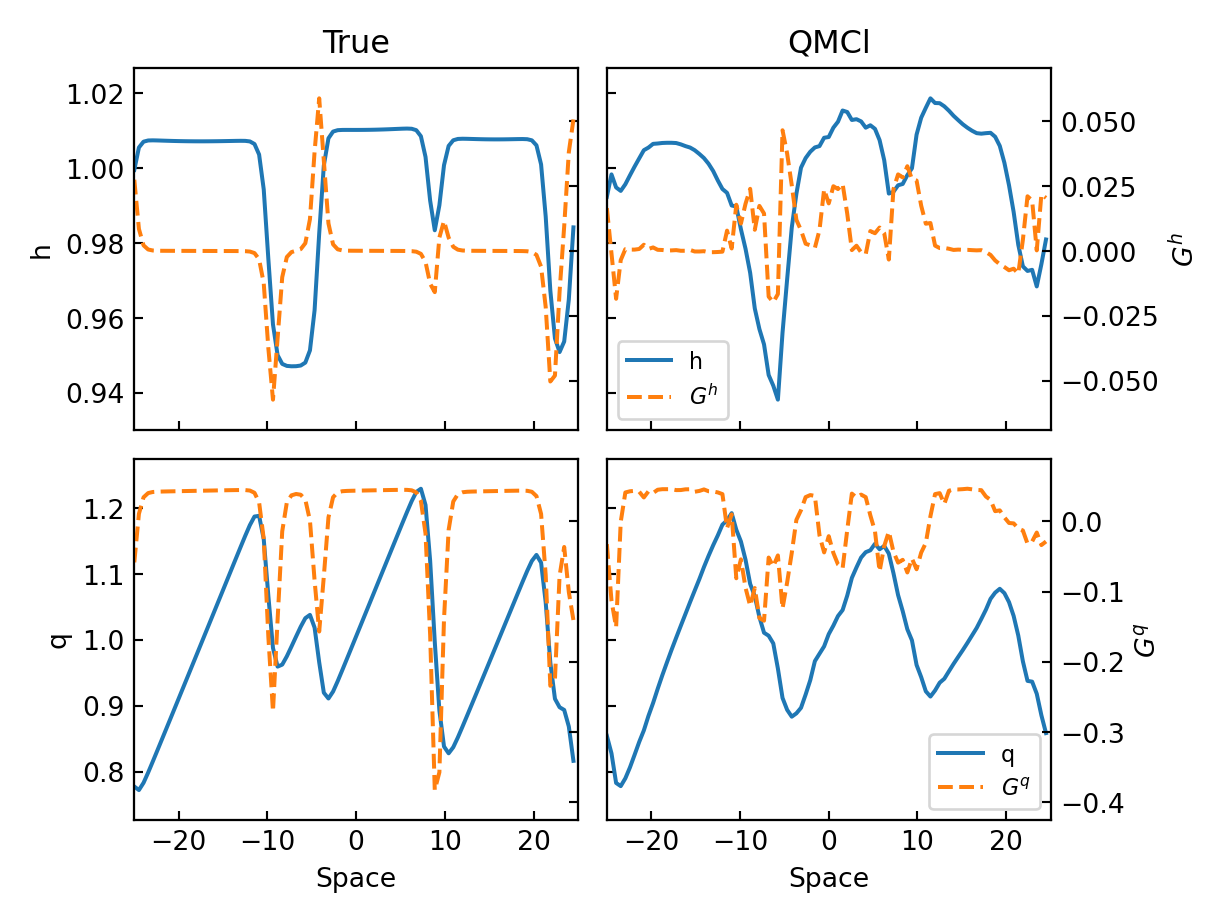}
    \caption{Spatial profiles of the final solution snapshot
    for the trajectory with $\delta=0.25$.
    The top row shows the results for the resolved state and the
    subgrid fluxes for the height field $h$;
    the bottom row for the momentum field $q$.
    In each row, the true results are on the left subplot,
    the predicted ones on the right subplot.
    Within each subplot, full lines are used for the resolved state,
    with their numerical values indicated on the left axis;
    dashed lines are used for the subgrid fluxes, with their
    values indicated on the right axis.}
    \label{fig:res01b}
\end{figure}

\begin{figure}[t]
    \centering
    \includegraphics[width=.7\textwidth]{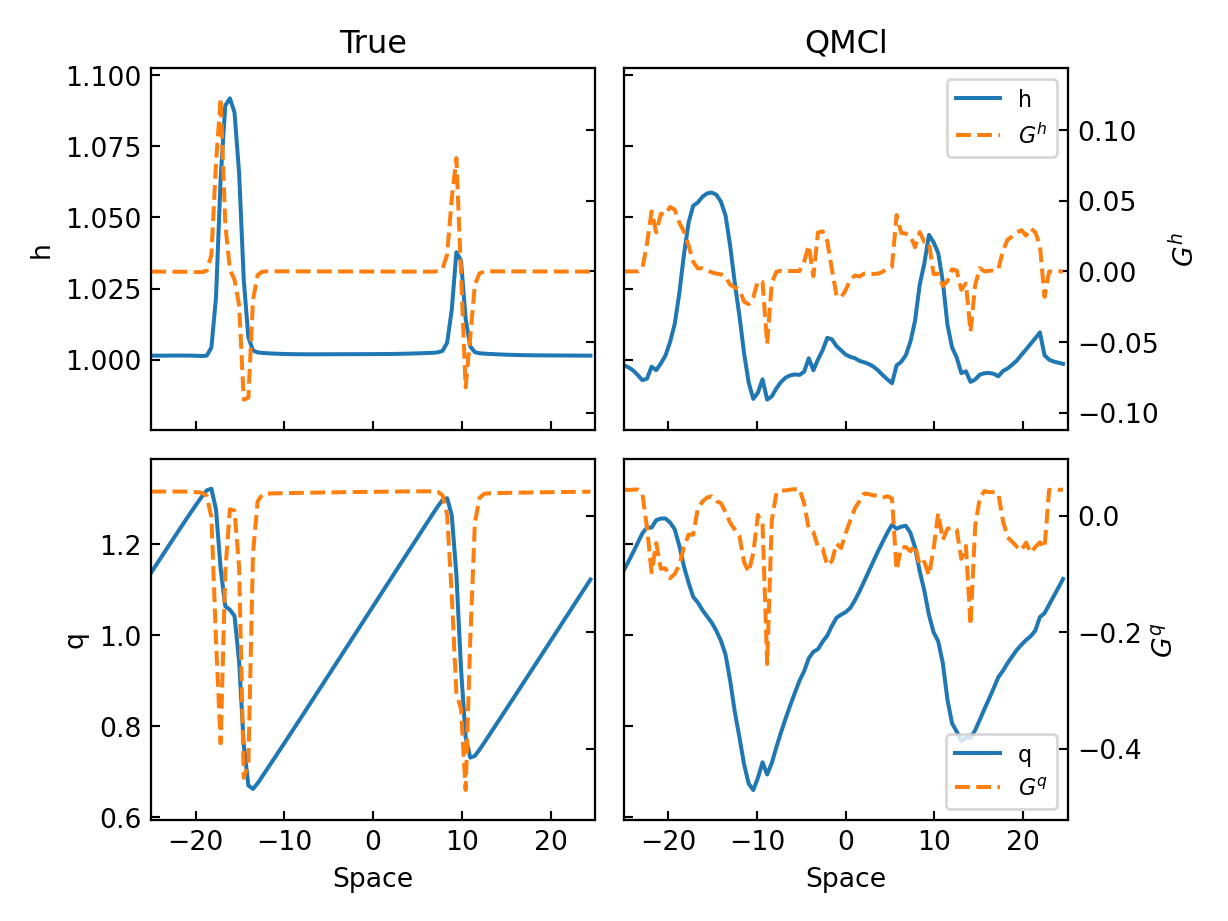}
    \caption{Spatial profiles of the final solution snapshot
    for the trajectory with $\delta=0.75$.
    The layout of the plots is as in Figure \ref{fig:res01b}.}
    \label{fig:res02b}
\end{figure}

The results comparing the evolution under the true and surrogate
models are grouped into four figures.
Figures \ref{fig:res01a} and \ref{fig:res01b}
include the results for the trajectory generated by the initial
condition obtained from \eqref{eq:ic-family} for $\delta=0.25$.
Figures \ref{fig:res02a} and \ref{fig:res02b}
showcase the results obtained for $\delta=0.75$.

\begin{figure}[t]
    \centering
    \includegraphics[width=.7\textwidth]{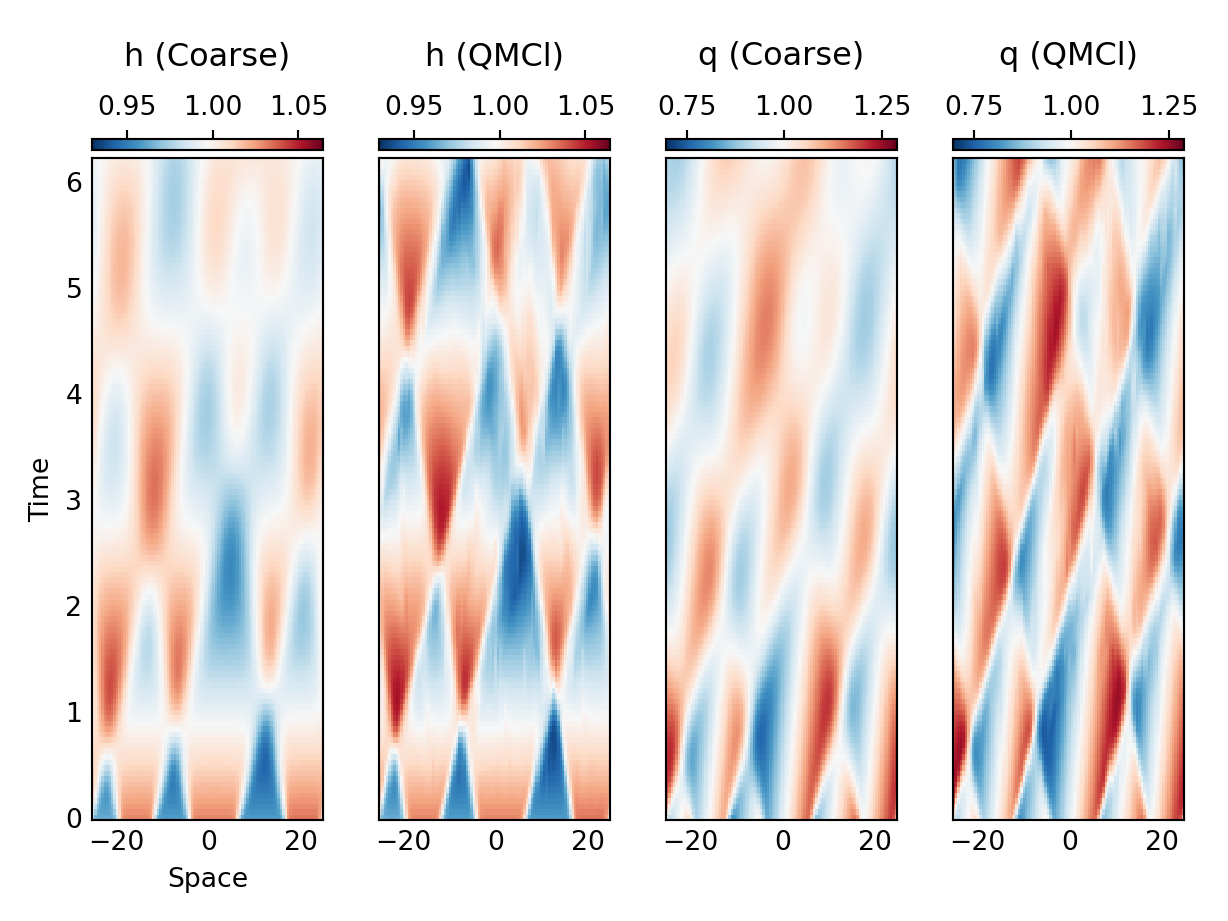}
    \caption{Space-time heatmaps for the trajectory with $\delta=0.25$
    comparing the dynamics on the coarse grid with the surrogate
    dynamics obtained by QMCl.
    From left to right, the first two plots show the results
    for the height field $h$, with the coarse grid results on the
    left and the QMCl results on the right;
    the next two plots show the results for the momentum field $q$,
    with coarse grid results on the left and QMCl on the right.
    The colorbars for each subplot are as in Figure \ref{fig:res01a}.}
    \label{fig:res01c}
\end{figure}

Looking at the heatmaps for the height and momentum fields in
Figure \ref{fig:res01a} for the trajectory with $\delta=0.25$
and in Figure \ref{fig:res02a} for the one with $\delta=0.75$,
we see that QMCl clearly captures the main features of the
dynamics for the tested time horizon.
For both trajectories, QMCl correctly predicts the main
traveling wave features of the dynamics, including the spatial
variations of the waves as they interact with one another.
Moreover, the magnitude of the predicted fields is consistent
with the true dynamics in most regions of the spatial domain.
In addition to the results shown in these figures, we have
also verified that QMCl can reproduce the dynamical trajectories
included in its training dataset with good accuracy.

In the present closure model, the subgrid fluxes play the role of
antidiffusive corrections to the fluxes generated entirely on the
coarse grid \cite{Timofeyev2025}.
Their values are proportional to the spatial derivative of the associated
field, taking large values when the height or momentum field is changing
rapidly in space, and approximately zero value when the fields are
constant in space.
As a result, the contribution of the subgrid fluxes is critical in
regions of strong spatial variation of the height and momentum fields,
counterbalancing the increased diffusion of the solution introduced by
our coarsening of the original dynamics by spatial averaging.
The strong spatial variation of the true subgrid fluxes renders the
considered dynamics challenging to resolve and predict accurately
by any closure scheme
(Section \ref{sec:swe-discussion}).
Typically, when the exact value of the subgrid fluxes is not
predicted correctly by the employed closure model (derived by
QMCl or another method), the surrogate evolution will be more
diffusive than the true evolution.

The effect of diffusion can be seen in the heatmaps for the
two trajectories.
By comparing the true and predicted fields of height and momentum,
we see that the surrogate dynamics involves stronger diffusion of
the solution in areas where there is strong spatial variation.
The increased diffusion is seen in the plots as bleeding of color
between different space-time regions, with the boundaries between
those regions being less sharp than in the true dynamics.
Similarly, the same effect can be seen by comparing the subgrid flux
fields for height and momentum.
In the true dynamics, the peaks of the fluxes are sharp and clearly
separated from the valleys, with clearly visible waves traveling in
opposite directions.
In the predicted dynamics, although the waves are still predicted
correctly, their peaks are not as sharp, with their values
diffusing more strongly to the surrounding region.

The subgrid flux plots also clearly demonstrate the effect of the
bayesian conditioning of the quantum density operators, which for
these results is performed every $10$ timesteps
($\approx 0.5$ time units).
One can see the discontinuities that occur every $10$ timesteps
for the predicted fluxes over large regions of the spatial domain.
These discontinuities represent corrections introduced to the predicted
flux values every time conditioning of the density operators is performed.
We perform the bayesian conditioning every $10$ timesteps to reduce
the computational cost of the simulation, since the conditioning is one
of the most expensive operations in our closure scheme
(Section \ref{sec:swe-discussion}).
In general, the prediction accuracy increases with increased
conditioning frequency.
Based on our numerical experiments, conditioning the density operators
every $10$ timesteps strikes a good balance between prediction accuracy
and computational cost for the considered dynamics.

Similar conclusions can be drawn from the spatial profiles shown
in Figure \ref{fig:res01b} for the trajectory with $\delta=0.25$
and Figure \ref{fig:res02b} for that with $\delta=0.75$.
In these plots one can see that the subgrid fluxes are proportional
to the spatial derivative of the associated field, leading to sharp
variations in space.
The predicted dynamics preserves the main features of the true
fields, but struggles to capture the sharp peaks of the subgrid
fluxes, leading to stronger diffusion of the height and momentum
fields in regions of strong spatial variation.

To further demonstrate the crucial role played by the subgrid
fluxes in counterbalancing diffusion,
in Figure \ref{fig:res01c}
we compare the evolution of the QMCl surrogate dynamics
with the true dynamics generated entirely on the coarse grid,
namely with the subgrid fluxes \eqref{eq:subgrid-flux}
set to zero.
The comparison is carried out for the trajectory generated from
the initial condition \eqref{eq:ic-family} with $\delta=0.25$
(Figures \ref{fig:res01a} and \ref{fig:res01b}),
but similar conclusions can be drawn from other initial conditions.
The comparison demonstrates that the surrogate fluxes predicted
by QMCl successfully counteract the diffusive effect of the
spatial coarsening seen in the coarse grid simulation
for the tested time horizon.
Namely, although the subgrid fluxes predicted by QMCl are not as
accurate as the true subgrid fluxes in regions of strong spatial
variation of the solution, they are nevertheless effective
in counterbalancing diffusion and thus improving the accuracy
of the predicted dynamics.

\subsection{Discussion}\label{sec:swe-discussion}
In the following paragraphs we discuss several aspects of
the presented closure scheme and how they inform future
research directions.

\subsubsection*{Computational aspects}
Extending the QMCl framework to spatiotemporal dynamics comes with
the challenge of increased computational cost for its implementation,
primarily due to the increased size of the training dataset required
to capture the main features of dynamics that evolves in both space
and time.
One key area of focus for future work is reducing the cost
of the numerical implementation of QMCl and improving its scaling
with the number of temporal samples $N$, grid size $M$
and spectral resolution $L$.

In the offline stage, the most expensive operation is computing the
leading $L$ eigenfunctions of the $NM\times NM$ symmetric kernel
matrix $\matr{K}_N$ representing the kernel integral operator
\eqref{eq:data-integral-op}.
Even for moderate size parameters $N$ and $M$, solving such an
eigenvalue problem can quickly become very expensive.
As explained in Section \ref{sec:swe-numerics}, to address this
we performed a partial Cholesky factorization of $\matr{K}_N$
to approximate it by a lower rank matrix and solve the resulting
eigenvalue problem with reduced cost
\cite{YChen2024,Vales2025evd}.

As an alternative we also tried using the random Fourier
features (RFF) method to approximate $\matr{K}_N$
by a lower rank matrix \cite{Rahimi2007}.
However, we found that a relatively large number of Fourier
features was required to obtain a satisfactory approximation
of $\matr{K}_N$, compared to the rank parameter $r$
required for the partial Cholesky factorization
\cite{Yang2012}.
The poor approximations obtained using the RFF method led to
issues in the downstream bistochastic normalization of
the approximate kernel matrix, which is why we did not use this
approach in the present work.

In the online stage, the most expensive operation is the bayesian
conditioning of the quantum density operators; more specifically,
building the quantum effect operator from the current resolved state.
To form the quantum effect operators we use kernel
\eqref{eq:data-kernel-cond}
to compare the current resolved product state vector to all product
state vectors present in our training dataset.
This is the only online operation that scales with the size
of the training dataset $NM$, and requires that we have the
training data loaded in memory throughout the online stage.

One way of alleviating this cost is to form a reduced representation
of the training dataset and use it to build the quantum effect
operators.
For instance, one can employ the sampling algorithm underlying the
partial Cholesky factorization algorithms \cite{YChen2024}
to extract a representative subset of the original training
dataset, and rely only on that subset for the conditioning of
the density operators.
An alternative would be to use the RFF method to approximate the
out of sample evaluation of kernel \eqref{eq:data-kernel-cond},
thereby eliminating the dependence on the training data
\cite{Rahimi2007,Giannakis2023}.
In our case, this would require keeping in memory one
$L\times L$ matrix for each random Fourier feature used,
leading to an unsustainable memory cost.

In QMCl the surrogate fluxes are computed as the expected value
of an observable with respect to a density operator.
As a result, to compute fluxes that vary strongly over the
spatial domain one has to ensure that the used density operators
are very \emph{sharp}: obtaining a clear maximum in some
regions of product space $\Omega$ and being close to zero
elsewhere, formally similar to a sum of Dirac delta measures.
Working with such sharp density operators poses a number of
computational challenges.

Producing accurate but very sharp density operators requires
a very precise choice of the bandwidth parameter $\epsilon'$
used in kernel \eqref{eq:data-kernel-cond}
for the bayesian conditioning.
Having such a selective kernel means that the resulting density
operators lead to highly accurate fluxes when the correct part of
the training dataset is identified as dynamically relevant.
However, it can lead to inaccurate fluxes when regions of the
training dataset are misidentified as relevant, compared to using a
less selective kernel that leads to wider averaging of the employed
observable.
Moreover, working with a very selective kernel makes it more likely
that no part of the training dataset is identified as dynamically
relevant, leading to an approximately zero feature vector and
posterior density operator.
In such a case, attempting to normalize the posterior density operator
to unit norm can yield unpredictable results.
Working with a wider, less selective kernel makes this scenario less
likely to occur, albeit with the tradeoff of less informative
conditioning and thus less accurate predicted fluxes.
Finally, even if very sharp density operators can be consistently
obtained, their accurate representation with respect to our eigenfunction
basis will likely require a large number of eigenfunctions, leading to
more expensive simulations.

The challenges involved in producing and working with very
selective kernel functions motivate our use of a kernel
bandwidth tuning algorithm \cite{Coifman2008}
and of the introduction of variable bandwidth
\cite{Berry2016,Giannakis2019acha}.
These measures help us define a kernel function that is
tailored to the characteristics of our training dataset,
which in turn allows us to strike a reasonable
compromise between numerical stability and sharpness
of the kernel used for bayesian conditioning.
As a result, we are able to construct a closure scheme
that can predict the main features of the considered
dynamics with reasonable computational cost and without
numerical instabilities.

\subsubsection*{Theoretical aspects}
One of the distinguishing features of the QMCl framework is that
the classical dynamics is first embedded into the noncommutative
setting of operator algebras used in quantum mechanics,
and only then discretized to finite dimension.
Namely, instead of directly discretizing functions on
Hilbert space $H$ to obtain functions on finite dimensional
subspace $H^L$,
we first map functions on $H$ to operators in $B(H)$
and then discretize them to operators in $B(H^L)$.
This process renders the discretization employed by QMCl
positivity preserving, with positive classical observables
and probability density functions mapped to positive definite
operators.

In addition to positivity preservation, this discretization
process embeds the discretized quantities in a space of larger
dimension than direct discretization.
More specifically, direct discretization maps real functions
on $H$ to functions on the $L$-dimensional subspace $H^L$.
The QMCl discretization maps real functions on $H$ to selfadjoint
real operators in $B(H^L)$, which is a subspace of dimension
$\frac{1}{2}L(L+1)$.
A potential line of future research is exploring whether this
difference in dimension, along with the noncommutative nature
of $B(H^L)$, allow QMCl to encode more statistical information
about the original dynamics compared to closure schemes performed
directly in $H^L$.

Another line of research concerns the out of sample performance of QMCl.
In Sections \ref{sec:qmcl-multi} and \ref{sec:swe-numerics}
we considered training data generated from a finite collection of
initial conditions, without making reference to an explicit process
of generating these initial conditions.
A natural next step is to consider an invariant measure $\mu$ that
admits a decomposition into conditional measures $\mu_i$,
sampled from a probability distribution of initial conditions.
Interesting questions in this direction are to characterize
the behavior of our closure scheme as the number of sampled
trajectories increases, and to derive associated out of sample
error bounds.

\section{Conclusion}\label{sec:conclusion}
We developed a framework for the dynamical closure of
spatiotemporal dynamics which combines elements of
the quantum mechanical closure (QMCl) framework developed in
\cite{Freeman2024}
and the vector valued spectral analysis (VSA) methodology of
\cite{Giannakis2019vsa}.
The presented framework encodes statistical information about the
unresolved degrees of freedom of the considered dynamics in the form
of a field of quantum density operators on the spatial domain.
This results in a flexible architecture where density operators are
collocated with the discretization cells of the coarse model
for the resolved dynamics.
The density operators communicate with the resolved dynamics by:
(1) providing the required flux terms through the expected value of
quantum observables;
(2) receiving information through kernel feature maps.
Importantly, our approach inherits the symmetry factorization properties
of VSA during the basis construction process, and can also be used with
training data consisting of multiple dynamical trajectories.
As with the original QMCl work, the process of embedding the
original dynamics into the noncommutative setting of quantum mechanics
leads to the positivity preserving discretization of classical
observables and probability density functions.

We applied QMCl to a closure problem for the shallow water
equations on a periodic one dimensional domain.
The numerical results demonstrate that QMCl can effectively
predict the dynamics of the out of sample initial conditions
tested in this work, generating height and momentum fields that
agree both qualitatively and quantitatively with the true dynamics.
As expected from the strong spatial variation of the subgrid fluxes,
the QMCl surrogate dynamics is more diffusive than the true dynamics,
and eventually struggles to fully capture the sharp spatial features
of the true solution.
Nevertheless, we believe that the results act as a successful proof
of concept for the extension of QMCl to spatiotemporal dynamics,
especially given the modest amount of data used for its training.
Everything else kept the same, we expect that increasing the
number of training samples or the spectral resolution will lead
to improved prediction accuracy.

Future research directions include improving the computational cost
of QMCl and the development of strategies for sampling initial conditions
in the training stage.
In addition, we plan to explore applications of QMCl to systems
exhibiting chaotic spatiotemporal dynamics.

\section*{Acknowledgements}
DG and JS acknowledge support from the U.S.
Department of Energy under grant DE-SC0025101.
CV was supported as a postdoctoral researcher from this grant.
DG acknowledges support from the U.S. Department of Defense,
Basic Research Office under Vannevar Bush Faculty Fellowship
grant N00014-21-1-2946.
DCF was supported as a PhD student from this grant.
We thank Ilon Joseph, Michael Montgomery and Ilya Timofeyev
for helpful conversations on the topics of this work.

\appendix
\section{Numerical implementation}\label{app:implementation}
The process of building and using the developed QMCl scheme to simulate
the surrogate dynamics consists of two stages:
the offline and online stages.
Below we outline the steps followed in each stage.

\subsubsection*{Offline stage}
The offline stage consists of the operations needed to build the
QMCl scheme using the training samples for the resolved state variables
$\{\hat{u}_n\}_{n=0}^{N-1}$
and the subgrid fluxes
$\{f(u_n,\cdot)\}_{n=0}^{N-1}$.
This involves the following sequence of steps.
\begin{enumerate}
\item Use the state samples $\{\hat{u}_n\}$
to calibrate the bandwidth parameter $\epsilon$
for the kernel \eqref{eq:data-kernel}
and perform its bistochastic normalization.
\item Use the state samples $\{\hat{u}_n\}$
to compute the leading $L$ eigenfunctions $\{\phi_\ell\}$
of the normalized integral operator
\eqref{eq:data-integral-op},
defining the subspace $H_N^L$ as their linear span.
The eigenfunction basis is represented by the $NM\times L$
matrix $\matr{\Phi}$ holding one eigenfunction in each column.
\item Form the transfer operator $P_N^L$ and its
$L\times L$ coefficients matrix $\matr{P}_N^L$.
\item Use the flux samples $\{f(u_n,\cdot)\}$
to form the two quantum observables
$A_0$, $A_1\in B(H_N^L)$,
computing their $L\times L$ coefficients matrices
$\matr{A}_0$ and $\matr{A}_1$.
Observable $A_0$ is built using the flux samples
$\{f_0(x_n,\cdot)\}$
corresponding to the height $h$ variables;
observable $A_1$ using $\{f_1(x_n,\cdot)\}$
corresponding to the $q$ variables.
\item Define the field of quantum density operators
$\rho\colon S\to B(H_N^L)$,
assigning one density operator to each coarse grid cell.
For each gridpoint $x_m\in S_M$, the associated density operator
$\rho(x_m)$ is represented by the vector
$\vect{\rho}_m$ in $\Rbb^L$.
\item Use the state samples $\{\hat{u}_n\}$
to calibrate the bandwidth parameter $\epsilon'$
of the kernel \eqref{eq:data-kernel-cond}
and compute the bandwidth function needed for its variable
bandwidth generalization.
\end{enumerate}

\subsubsection*{Online stage}
In the online stage we use the model built in the offline stage
to timestep the surrogate dynamics.
We denote by $\tilde{u}\in X_r$ the resolved state predicted by
QMCl, to distinguish it from the true resolved state $\hat{u}$.
The online stage consists of the following steps.
\begin{enumerate}
\item Initialize the resolved state $\tilde{u}$.
\item Initialize the density operators $\rho$ representing
    each coarse grid cell by using an uninformative initial condition:
    set each density operator $\rho(x_m)$ to the projection along the
    unit function $1_N\in H_N^L$ that takes the value $1$ everywhere.
\item Condition each density operator using the initial resolved state
    $\tilde{u}$ and the training data $\{\hat{u}_n\}$.
\item Use the density operators $\rho$ and the observables
    $A_0$ and $A_1$ to compute the surrogate fluxes $\tilde{f}$,
    corresponding to two flux values for each coarse grid cell.
\item Begin the timestepping loop.
\begin{enumerate}
\item Update the resolved state $\tilde{u}$ using the surrogate fluxes.
\item Update the density operators $\rho$ using the transfer operator.
\item Condition the density operators $\rho$ using the updated resolved
    state $\tilde{u}$ and the state training samples $\{\hat{u}_n\}$.
\item Compute the updated surrogate fluxes $\tilde{f}$.
\end{enumerate}
\end{enumerate}

In the timestepping loop as presented above, the bayesian
conditioning of the density operators is performed every
timestep.
In our presented numerical results, we have instead elected
to condition the density operators every $10$ timesteps
(Section \ref{sec:swe-numerics}).

\bibliographystyle{elsarticle-num}
\biboptions{sort&compress}
\bibliography{refs}

\end{document}